\documentclass[reqno,11pt]{amsart}
\usepackage{amsmath}
\usepackage{amsxtra}
\usepackage{amscd}
\usepackage{amsthm}
\usepackage{amsfonts}
\usepackage{amssymb}
\usepackage{eucal}
\usepackage[matrix,arrow,curve]{xy}

\theoremstyle{plain}
\newtheorem{Thm}[subsection]{Theorem}
\newtheorem{Cor}[subsection]{Corollary}
\newtheorem{Lem}[subsection]{Lemma}
\newtheorem{Prop}[subsection]{Proposition}
\newtheorem{Conj}[subsection]{Conjecture}

\theoremstyle{definition}
\newtheorem{Def}[subsection]{Definition}

\theoremstyle{remark}

\newtheorem{Rem}[subsection]{Remark}

\newtheorem{conj}{Conjecture}

\numberwithin{equation}{section}
\renewcommand{\rm}{\normalshape}

\newif\ifShowLabels
\ShowLabelstrue
\newdimen\theight
\def\TeXref#1{%
    \leavevmode\vadjust{\setbox0=\hbox{{\tt
        \quad\quad  {\small \rm #1}}}%
    \theight=\ht0
    \advance\theight by \lineskip
    \kern -\theight \vbox to
    \theight{\rightline{\rlap{\box0}}%
    \vss}%
    }}%

\ShowLabelsfalse

\renewcommand{\sec}[2]{\section{#2}\label{S:#1}%
    \ifShowLabels \TeXref{{S:#1}} \fi}
\newcommand{\ssec}[2]{\subsection{#2}\label{SS:#1}%
    \ifShowLabels \TeXref{{SS:#1}} \fi}

\newcommand{\sssec}[2]{\subsubsection{#2}\label{SSS:#1}%
    \ifShowLabels \TeXref{{SSS:#1}} \fi}

\newcommand{\refs}[1]{Section ~\ref{S:#1}}
\newcommand{\refss}[1]{Section ~\ref{SS:#1}}
\newcommand{\refsss}[1]{Section ~\ref{SSS:#1}}
\newcommand{\reft}[1]{Theorem ~\ref{T:#1}}
\newcommand{\refl}[1]{Lemma ~\ref{L:#1}}
\newcommand{\refp}[1]{Proposition ~\ref{P:#1}}

\newcommand{\refd}[1]{Definition ~\ref{D:#1}}

\newcommand{\refe}[1]{\eqref{E:#1}}
\newcommand{\refco}[1]{Conjecture ~\ref{Co:#1}}

\newenvironment{thm}[1]%
    { \begin{Thm} \label{T:#1}  \ifShowLabels \TeXref{T:#1} \fi }%
    { \end{Thm} }

\renewcommand{\th}[1]{\begin{thm}{#1} \sl }
\renewcommand{\eth}{\end{thm} }

\newenvironment{lemma}[1]%
    { \begin{Lem} \label{L:#1}  \ifShowLabels \TeXref{L:#1} \fi }%
    { \end{Lem} }
\newcommand{\lem}[1]{\begin{lemma}{#1} \sl}
\newcommand{\elem}{\end{lemma}}

\newenvironment{propos}[1]%
    { \begin{Prop} \label{P:#1}  \ifShowLabels \TeXref{P:#1} \fi }%
    { \end{Prop} }
\newcommand{\prop}[1]{\begin{propos}{#1}\sl }
\newcommand{\eprop}{\end{propos}}

\newenvironment{corol}[1]%
    { \begin{Cor} \label{C:#1}  \ifShowLabels \TeXref{C:#1} \fi }%
    { \end{Cor} }
\newcommand{\cor}[1]{\begin{corol}{#1} \sl }
\newcommand{\ecor}{\end{corol}}

\newenvironment{defeni}[1]%
    { \begin{Def} \label{D:#1}  \ifShowLabels \TeXref{D:#1} \fi }%
    { \end{Def} }
\newcommand{\defe}[1]{\begin{defeni}{#1} \sl }
\newcommand{\edefe}{\end{defeni}}

\newenvironment{remark}[1]%
    { \begin{Rem} \label{R:#1}  \ifShowLabels \TeXref{R:#1} \fi }%
    { \end{Rem} }
\newcommand{\rem}[1]{\begin{remark}{#1}}
\newcommand{\erem}{\end{remark}}

\newenvironment{conjec}[1]%
    { \begin{Conj} \label{Co:#1}  \ifShowLabels \TeXref{Co:#1} \fi }%
    { \end{Conj} }
\renewcommand{\conj}[1]{\begin{conjec}{#1} \sl }
\newcommand{\econj}{\end{conjec}}

\newcommand{\eq}[1]%
    { \ifShowLabels \TeXref{E:#1} \fi
       \begin{equation} \label{E:#1} }
\newcommand{\eeq}{ \end{equation} }

\newcommand{\prf}{ \begin{proof} }
\newcommand{\epr}{ \end{proof} }

\newcommand{\nc}{\newcommand}

\newcommand{\iso}{\stackrel{\sim}{\longrightarrow}}
\newcommand{\lra}{\longrightarrow}
\newcommand{\hra}{\hookrightarrow}

\nc{\HC}{{\mathcal{HC}}}
\nc{\on}{\operatorname}
\nc{\BA}{{\mathbb{A}}}
\nc{\BC}{{\mathbb{C}}}
\nc{\BG}{{\mathbb{G}}}
\nc{\BM}{{\mathbb{M}}}
\nc{\BN}{{\mathbb{N}}}
\nc{\BQ}{{\mathbb{Q}}}
\nc{\BP}{{\mathbb{P}}}
\nc{\BR}{{\mathbb{R}}}
\nc{\BZ}{{\mathbb{Z}}}
\nc{\BS}{{\mathbb{S}}}

\nc{\CA}{{\mathcal{A}}}
\nc{\CB}{{\mathcal{B}}}
\nc{\CalC}{{\mathcal C}}
\nc{\CalD}{{\mathcal D}}
\nc{\CE}{{\mathcal{E}}}
\nc{\CF}{{\mathcal{F}}}
\nc{\CG}{{\mathcal{G}}}
\nc{\CH}{{\mathcal{H}}}
\nc{\CI}{{\mathcal{I}}}
\nc{\CJ}{{\mathcal{J}}}
\nc{\CK}{{\mathcal{K}}}
\nc{\CL}{{\mathcal{L}}}
\nc{\CM}{{\mathcal{M}}}
\nc{\CMM}{{\mathcal{M}^{\operatorname{gen}}_\hbar(-\rho)}}
\nc{\CN}{{\mathcal{N}}}
\nc{\CO}{{\mathcal{O}}}
\nc{\CP}{{\mathcal{P}}}
\nc{\CQ}{{\mathcal{Q}}}
\nc{\CR}{{\mathcal{R}}}
\nc{\CS}{{\mathcal{S}}}
\nc{\CT}{{\mathcal{T}}}
\nc{\CU}{{\mathcal{U}}}
\nc{\CV}{{\mathcal{V}}}
\nc{\CW}{{\mathcal{W}}}
\nc{\CX}{{\mathcal{X}}}
\nc{\CY}{{\mathcal{Y}}}
\nc{\CZ}{{\mathcal{Z}}}

\nc{\gen}{{\operatorname{gen}}}
\nc{\cM}{{\check{\mathcal M}}{}}
\nc{\csM}{{\check{\mathcal A}}{}}
\nc{\obM}{{\overset{\circ}{\mathbf M}}{}}
\nc{\oCA}{{\overset{\circ}{\mathcal A}}{}}
\nc{\obA}{{\overset{\circ}{\mathbf A}}{}}
\nc{\ooM}{{\overset{\circ}{M}}{}}
\nc{\osM}{{\overset{\circ}{\mathsf M}}{}}
\nc{\vM}{{\overset{\bullet}{\mathcal M}}{}}
\nc{\nM}{{\underset{\bullet}{\mathcal M}}{}}
\nc{\obD}{{\overset{\circ}{\mathbf D}}{}}
\nc{\cp}{{\overset{\circ}{\mathbf p}}{}}
\nc{\ofZ}{{\overset{\circ}{\mathfrak Z}}{}}

\nc{\fa}{{\mathfrak{a}}}
\nc{\fb}{{\mathfrak{b}}}
\nc{\fg}{{\mathfrak{g}}}
\nc{\fgl}{{\mathfrak{gl}}}
\nc{\fh}{{\mathfrak{h}}}
\nc{\fj}{{\mathfrak{j}}}
\nc{\fl}{{\mathfrak{l}}}
\nc{\fm}{{\mathfrak{m}}}
\nc{\fn}{{\mathfrak{n}}}
\nc{\fu}{{\mathfrak{u}}}
\nc{\fp}{{\mathfrak{p}}}
\nc{\frr}{{\mathfrak{r}}}
\nc{\fs}{{\mathfrak{s}}}
\nc{\ft}{{\mathfrak{t}}}
\nc{\fT}{{\mathfrak{T}}}
\nc{\ofT}{{\overline{\mathfrak T}}}
\nc{\ofS}{{\overline{\mathfrak S}}}
\nc{\fsl}{{\mathfrak{sl}}}
\nc{\hsl}{{\widehat{\mathfrak{sl}}}}
\nc{\hgl}{{\widehat{\mathfrak{gl}}}}
\nc{\hg}{{\widehat{\mathfrak{g}}}}
\nc{\chg}{{\widehat{\mathfrak{g}}}{}^\vee}
\nc{\hn}{{\widehat{\mathfrak{n}}}}
\nc{\chn}{{\widehat{\mathfrak{n}}}{}^\vee}

\nc{\fA}{{\mathfrak{A}}}
\nc{\fB}{{\mathfrak{B}}}
\nc{\fD}{{\mathfrak{D}}}
\nc{\fE}{{\mathfrak{E}}}
\nc{\fF}{{\mathfrak{F}}}
\nc{\fG}{{\mathfrak{G}}}
\nc{\fI}{{\mathfrak{I}}}
\nc{\fJ}{{\mathfrak{J}}}
\nc{\fK}{{\mathfrak{K}}}
\nc{\fL}{{\mathfrak{L}}}
\nc{\fM}{{\mathfrak{M}}}
\nc{\fN}{{\mathfrak{N}}}
\nc{\frP}{{\mathfrak{P}}}
\nc{\fQ}{{\mathfrak Q}}
\nc{\fS}{{\mathfrak S}}
\nc{\fU}{{\mathfrak{U}}}
\nc{\fV}{{\mathfrak{V}}}
\nc{\fZ}{{\mathfrak{Z}}}

\nc{\bb}{{\mathbf{b}}}
\nc{\bc}{{\mathbf{c}}}
\nc{\be}{{\mathbf{e}}}
\nc{\bj}{{\mathbf{j}}}
\nc{\bn}{{\mathbf{n}}}
\nc{\bp}{{\mathbf{p}}}
\nc{\bq}{{\mathbf{q}}}
\nc{\bs}{{\mathbf{s}}}
\nc{\bt}{{\mathbf{t}}}
\nc{\bu}{{\mathbf{u}}}
\nc{\bv}{{\mathbf{v}}}
\nc{\bx}{{\mathbf{x}}}
\nc{\by}{{\mathbf{y}}}
\nc{\bw}{{\mathbf{w}}}
\nc{\bA}{{\mathbf{A}}}
\nc{\bB}{{\mathbf{B}}}
\nc{\bC}{{\mathbf{C}}}
\nc{\bK}{{\mathbf{K}}}
\nc{\bD}{{\mathbf{D}}}
\nc{\bH}{{\mathbf{H}}}
\nc{\bI}{{\mathbf{I}}}
\nc{\bM}{{\mathbf{M}}}
\nc{\bN}{{\mathbf{N}}}
\nc{\bO}{{\mathbf{O}}}
\nc{\bQ}{{\mathbf Q}}
\nc{\bS}{{\mathbf{S}}}
\nc{\bT}{{\mathbf{T}}}
\nc{\bU}{{\mathbf{U}}}
\nc{\bV}{{\mathbf{V}}}
\nc{\bW}{{\mathbf{W}}}
\nc{\bX}{{\mathbf{X}}}
\nc{\bP}{{\mathbf{P}}}
\nc{\bZ}{{\mathbf{Z}}}
\nc{\bmu}{{\boldsymbol{\mu}}}

\nc{\sA}{{\mathsf{A}}}
\nc{\sB}{{\mathsf{B}}}
\nc{\sC}{{\mathsf{C}}}
\nc{\sD}{{\mathsf{D}}}
\nc{\sF}{{\mathsf{F}}}
\nc{\sK}{{\mathsf{K}}}
\nc{\sM}{{\mathsf{M}}}
\nc{\sO}{{\mathsf{O}}}
\nc{\sQ}{{\mathsf{Q}}}
\nc{\sP}{{\mathsf{P}}}
\nc{\sT}{{\mathsf{T}}}
\nc{\sV}{{\mathsf{V}}}
\nc{\sW}{{\mathsf{W}}}
\nc{\sZ}{{\mathsf{Z}}}
\nc{\sfp}{{\mathsf{p}}}
\nc{\sr}{{\mathsf{r}}}
\nc{\st}{{\mathsf{t}}}
\nc{\sv}{{\mathsf{v}}}
\nc{\sw}{{\mathsf{w}}}
\nc{\sx}{{\mathsf{x}}}
\nc{\sfb}{{\mathsf{b}}}
\nc{\sfc}{{\mathsf{c}}}
\nc{\sd}{{\mathsf{d}}}
\nc{\sg}{{\mathsf{g}}}
\nc{\sfl}{{\mathsf{l}}}

\nc{\BK}{{\bar{K}}}

\nc{\tA}{{\widetilde{\mathbf{A}}}}
\nc{\tB}{{\widetilde{\mathcal{B}}}}
\nc{\tg}{{\widetilde{\mathfrak{g}}}}
\nc{\tG}{{\widetilde{G}}}
\nc{\TM}{{\widetilde{\mathbb{M}}}{}}
\nc{\tO}{{\widetilde{\mathsf{O}}}{}}
\nc{\tU}{{\widetilde{\mathfrak{U}}}{}}
\nc{\TZ}{{\tilde{Z}}}
\nc{\tZ}{\widetilde{Z}{}}
\nc{\tx}{{\tilde{x}}}
\nc{\tbv}{{\tilde{\bv}}}
\nc{\tfP}{{\widetilde{\mathfrak{P}}}{}}
\nc{\tz}{{\tilde{\zeta}}}
\nc{\tmu}{{\tilde{\mu}}}

\nc{\td}{\ddot{\underline{d}}{}}
\nc{\tzeta}{\widetilde{\zeta}{}}
\nc{\hd}{{\widehat{\underline{d}}}}
\nc{\hG}{{\widehat{G}}}
\nc{\hBP}{\widehat{\mathbb P}{}}
\nc{\hQ}{{\widehat{Q}}}
\nc{\hsM}{\widehat{\mathsf M}{}}
\nc{\hfM}{\widehat{\mathfrak M}{}}
\nc{\hCP}{\widehat{\mathcal P}{}}
\nc{\hCR}{\widehat{\mathcal R}{}}
\nc{\hCS}{{\widehat{\mathcal S}}}
\nc{\hfZ}{\widehat{\mathfrak Z}{}}

\nc{\urho}{\underline{\rho}}
\nc{\uB}{\underline{B}}
\nc{\uC}{{\underline{\mathbb{C}}}}
\nc{\ui}{\underline{i}}
\nc{\ofP}{{\overline{\mathfrak{P}}}}

\nc{\hrho}{{\hat{\rho}}}

\nc{\unl}{\underline}
\nc{\ol}{\overline}
\nc{\one}{{\mathbf{1}}}
\nc{\two}{{\mathbf{t}}}

\nc{\Tot}{{\mathop{\operatorname{\rm Tot}}}}
\nc{\Hilb}{{\mathop{\operatorname{\rm Hilb}}}}
\nc{\End}{{\operatorname{End}}}
\nc{\CHom}{{\mathop{\operatorname{{\mathcal{H}}\it om}}}}
\nc{\defi}{{\mathop{\operatorname{\rm def}}}}
\nc{\length}{{\mathop{\operatorname{\rm length}}}}

\nc{\Cliff}{{\mathsf{Cliff}}}
\nc{\Fib}{{\mathsf{Fib}}}
\nc{\Coh}{{\mathsf{Coh}}}
\nc{\FCoh}{{\mathsf{FCoh}}}

\nc{\reg}{{\on{reg}}}
\nc{\Sp}{{\on{Sp}}}

\nc{\cplus}{{\mathbf{C}_+}}
\nc{\cminus}{{\mathbf{C}_-}}
\nc{\cthree}{{\mathbf{C}_*}}
\nc{\Qbar}{{\bar{Q}}}
\nc{\Fl}{{{\mathcal F}\ell}}
\nc{\bh}{{\bar{h}}}
\nc{\bOmega}{{\overline{\Omega}}}
\nc\tGr{\widetilde{\Gr}}
\nc{\ul}{\underline}
\nc{\seq}[1]{\stackrel{#1}{\sim}}
\nc\ogu{\overline{G/U}}
\nc\chlam{\check{\lam}}

\nc\St{\operatorname{St}}
\nc{\bLambda}{{\boldsymbol{\Lambda}}}
\nc\uS{\underline{S}}
\nc\QM{\mathcal{QM}}
\nc{\chmu}{\check{\mu}}
\nc{\CHH}{{\CH\!\!\CH}}
\author{Michael Finkelberg and Vadim Schechtman}
\title{Microlocal approach to Lusztig's symmetries}

\dedicatory{To Victor Ginzburg on his 60th birthday}

\begin{document}

\begin{abstract}
We reformulate the De Concini -- Toledano Laredo conjecture about the
monodromy of the Casimir connection in terms of a relation between Lusztig's
symmetries of quantum group modules and the monodromy in the vanishing cycles
of factorizable sheaves.
\end{abstract}

\maketitle 

\sec{Intro}{Introduction}

\ssec{noname}{} Let $\fg$ be a semisimple Lie algebra, and $\fh$ a Cartan
subalgebra. Let $\fh^\reg\subset\fh$ be the complement to the root hyperplanes
arrangement. For an integrable $\fg$-module $V$ C.~De~Concini and 
C.~Procesi~\cite{dcp} have introduced
an integrable {\em Casimir} connection with coefficients in the trivial vector
bundle $V\otimes\CO_{\fh^\reg}$ (it was later rediscovered by 
J.~Millson, V.~Toledano~Laredo~\cite{mt}
and J.~Felder, Y.~Markov, V.~Tarasov, A.~Varchenko~\cite{fmtv}) 
and conjectured that its monodromy can be expressed in terms of the action of 
the quantum Weyl group~\cite{lu},~\cite{so}
on the corresponding Weyl module $W_V$ over the corresponding quantum
group $\bU_v(\fg)$. This conjecture was later independently formulated
and proved by V.~Toledano~Laredo for $v$ in the formal neighbourhood of~1.
The key notion introduced in his proof was the notion of a  (quasi-){\it Coxeter
category}. The original definition of this notion is of combinatorial nature.
We suggest a more topological version of this definition in~\refs{Nota}.
It is a collection of local systems of restriction functors on the open
strata of hyperplane arrangements arising from the root hyperplanes of a root
system, compatible under Verdier specialization. This approach makes it clear
for example that the category of representations of a rational Cherednik algebra
carries a Coxeter structure, see~\cite{be}.

One of the main examples of a Coxeter category is a category of integrable representations
of a quantum group. We consider the category $\mathcal C$ of representations of Lusztig's
small quantum group. It has a geometric realization as the category ${\mathcal{FS}}$ of
factorizable sheaves~\cite{bfs}. This is the category of certain compatible collections of
perverse sheaves on the configuration spaces of a Riemann surface. One of our key observations
is that the category ${\mathcal{FS}}$ carries a natural Coxeter structure (in our topological
definition). 

We conjecture that the equivalence $\Phi\colon{\mathcal{FS}}\iso\mathcal C$ of~\cite{bfs}
takes the Coxeter structure on ${\mathcal{FS}}$ to the Coxeter structure on
$\mathcal C$. This is essentially a reformulation of
De~Concini -- Toledano~Laredo conjecture (hence it follows from the results
of V.~Toledano-Laredo for $v$ in the formal neighbourhood of $1$). 
Roughly, it says that the monodromy in the vanishing cycles of factorizable
sheaves acts by Lusztig's symmetries. 

\ssec{precise}{}
Let us formulate the last point of~\refss{noname} more precisely.
We choose a Borel subalgebra $\fh\subset\fb\subset\fg$.
The corresponding set of simple roots is denoted by $I$; 
for $i\in I$ the corresponding simple root is denoted $\alpha_i$.
We fix a Weyl group invariant symmetric bilinear form $?\cdot?$ on $\fh^*$
such that the square length of a {\em short} root is $\alpha_i\cdot\alpha_i=2$.

We fix a primitive root of unity $\zeta$ of degree $\ell$; for simplicity in
this introduction, we assume that $\ell$ is not divisible by 2 and 3.
We consider an integral dominant weight $\lambda\in\fh^*$ such that 
$0\leq\langle\lambda,\check\alpha_i\rangle<\ell$ for any $i\in I$ (pairings with
simple coroots).

For $\beta=\sum_{i\in I}b_i\alpha_i\in\BN[I]$ 
we consider the configuration space $\BA^\beta$ of colored
divisors on the complex affine line $\BA^1$. The open subspace
$\overset{\diamondsuit}\BA{}^\beta\subset\BA^\beta$ of multiplicity free 
divisors on $\BA^1\setminus\{0\}$ carries
a 1-dimensional local system $\CJ^\beta_\lambda$
with the following monodromies:
$\zeta^{-2\alpha_i\cdot\alpha_j}$ when a point of colour $i$ goes
counterclockwise around a point of colour $j\ne i$;
$-\zeta^{-\alpha_i\cdot\alpha_i}$ when two points of colour $i$ trade their
positions going around a halfcircle counterclockwise; and
$\zeta^{2\lambda\cdot\alpha_i}$ when a point of colour $i$ goes around 0
counterclockwise. We denote by
$\CI^\beta_\lambda$ the Goresky-MacPherson extension of $\CJ^\beta_\lambda$ to
$\BA^\beta$ (a perverse sheaf).

We have a pairing $\langle\cdot,\cdot\rangle\colon\fh\times\BA^\beta\to\BA^1$
given in the coordinates $(t_{i,r})_{i\in I}^{1\leq r\leq b_i}$ on $\BA^\beta$,
and $(z_j)_{j\in I}$ in the basis of fundamental coweights on $\fh$,
as follows: $\langle(z_j),(t_{i,r})\rangle=\sum_{i\in I}z_i\sum_{r=1}^{b_i}t_{i,r}$.
The vanishing cycles 
$\Phi^\beta_\lambda:=\Phi_{\langle\cdot,\cdot\rangle}\on{p}^\circ\CI^\beta_\lambda$
of the pullback of $\CI^\beta_\lambda$ to $\fh\times\BA^\beta$ is a perverse
sheaf supported on $\fh\simeq\fh\times\{\beta\cdot0\}$.

We conjecture that $\Phi^\beta_\lambda$
is smooth along $\fh^\reg\subset\fh$. In order to describe its monodromy on
$\fh^\reg$, recall that one of the main results of~\cite{bfs} is a canonical 
identification of the stalk $(\Phi^\beta_\lambda)_{C_0}$ at the fundamental Weyl
chamber in $\fh^\reg_\BR$ with the weight space $L^\lambda_{\lambda-\beta}$ of the
irreducible module with highest weight $\lambda$ over the Lusztig big
quantum group $\dot\bU_\zeta$ (note that the restriction of $L^\lambda$ to
the Lusztig small quantum group $\dot\bu_\zeta$ remains irreducible since
$\lambda$ is an $\ell$-restricted weight). We conjecture that the local system
$\Phi^\beta_\lambda|_{\fh^\reg}$ is given by the following representation of the 
fundamental groupoid of $\fh^\reg$: the stalk at a Weyl chamber $wC_0$ in
$\fh^\reg_\BR$ is $L^\lambda_{w(\lambda-\beta)}\ (w$ runs through the Weyl group $W$), 
and the half monodromies around
the walls are given by the Lusztig symmetries $T'_{i,\pm1}$ and $T''_{i,\pm1}$
(see~\refs{algeb} for details).

\ssec{contents}{}
Here is the outline of the paper.

In~\refs{Exa} we consider an elementary example of type $A_2$. We compare the action of
Lusztig's symmetries in the ``almost extremal'' weight spaces of integrable modules over the
corresponding quantum group (i.e.\ the weights obtained from the extremal ones by subtracting
a root) with the monodromy action in the vanishing cycles of related perverse sheaves on the
plane.

In~\refs{Nota} we propose a topological reformulation of Toledano Laredo's notion of Coxeter
category. It is very similar to Deligne's topological reformulation~\cite{dd} of the notion of
braided tensor category.

In~\refs{algeb} we recall the (algebraic) construction of the Coxeter structure on the category
of integrable modules over Lusztig's big quantum group, in terms of Lusztig's symmetries.
It gives rise to the Coxeter structure on the category $\mathcal C$ of representations of
Lusztig's small quantum group.

In~\refs{topolo} we recall very concisely the category of factorizable sheaves ${\mathcal{FS}}$
introduced
in~\cite{bfs} and the equivalence $\Phi\colon{\mathcal{FS}}\iso\mathcal C$
(see~\refss{rev1} and~\refss{rev2}). We also take an opportunity to correct a confusion
in~\cite{bfs} between Langlands dual types (see~\refss{err}). Then we go on to construct a
Coxeter structure on ${\mathcal{FS}}$. The construction goes through the De Rham realization
of ${\mathcal{FS}}$, and works only for $v$ sufficiently close to $1$, but we expect it to
work for arbitrary $v$. The construction also uses some
results on iterated specialization and microlocalization over hyperplane 
arrangements presented in~\refs{itera} which might be of independent interest.
A more systematic approach to these questions is developed in~\cite{fks}.

Finally, in~\refs{discu} we formulate the main conjecture that the equivalence
$\Phi\colon{\mathcal{FS}}\iso\mathcal C$ intertwines the Coxeter structures on 
${\mathcal{FS}}$ and on $\mathcal C$.

\ssec{ack}{Acknowledgments} We are grateful to R.~Fedorov, S.~Khoroshkin,
B.~Feigin, A.~Postnikov, G.~Rybnikov, L.~Rybnikov,
V.~Toledano Laredo, D.~Gaitsgory, D.~Kazhdan, M.~Kapranov, A.~Braverman,
R.~Bezrukavnikov, A.Varchenko for the inspiring discussions,
and M.~Kashiwara for an important reference.

In fact, this note arose
from a question asked by R.~Fedorov in the summer 2012.

The research of M.F. has been funded within the framework of the HSE University
Basic Research Program and the Russian Academic Excellence Project `5-100'.

\sec{Exa}{An example}

\ssec{alg}{Algebra}
We follow the notations of~\cite{lu}. Let $\bU$ be the quantum universal
enveloping algebra of type $A_2$, over the ring $\CA=\BZ[v^{\pm1}]$.
The positive (resp. negative) subalgebra $\bU^+$ (resp. $\bU^-$)
is generated by the divided powers $E_i^{(r)}$ (resp $F_i^{(r)}),\ i=1,2,\
r\in\BN$. Let $\Lambda=(\mu_1,\mu_2)\in\BN^2$ be a dominant highest weight
such that $\mu_1\geq1\leq\mu_2$, and $L(\Lambda)$ the corresponding integrable
$\bU$-module with the highest vector $\sv$.
We will be interested in the weight spaces
$L(\Lambda)_{(\mu_1-1,\mu_2-1)},\ L(\Lambda)_{(-\mu_1+1,\mu_1+\mu_2-2)},\
L(\Lambda)_{(\mu_1+\mu_2-2,-\mu_2+1)},\ L(\Lambda)_{(\mu_2-1,-\mu_1-\mu_2+2)},\\
L(\Lambda)_{(-\mu_1-\mu_2+2,\mu_1-1)},\ L(\Lambda)_{(-\mu_2+1,-\mu_1+1)}$ (these weights
form a single Weyl group orbit). They have canonical bases
$(F_1F_2\sv,F_2F_1\sv),\ (F_1^{(\mu_1)}F_2\sv,F_2F_1^{(\mu_1)}\sv),\
(F_2^{(\mu_2)}F_1\sv,F_1F_2^{(\mu_2)}\sv),\\
(F_2^{(\mu_1+\mu_2-1)}F_1^{(\mu_1)}\sv,F_1F_2^{(\mu_1+\mu_2-1)}F_1^{(\mu_1-1)}\sv),\
(F_1^{(\mu_1+\mu_2-1)}F_2^{(\mu_2)}\sv,F_2F_1^{(\mu_1+\mu_2-1)}F_2^{(\mu_2-1)}\sv),\\
(F_1^{(\mu_2-1)}F_2^{(\mu_1+\mu_2-1)}F_1^{(\mu_1)}\sv=
F_2^{(\mu_1)}F_1^{(\mu_1+\mu_2-1)}F_2^{(\mu_2-1)}\sv,\\
F_1^{(\mu_2)}F_2^{(\mu_1+\mu_2-1)}F_1^{(\mu_1-1)}\sv=
F_2^{(\mu_1-1)}F_1^{(\mu_1+\mu_2-1)}F_2^{(\mu_2)}\sv)$, respectively.

We are interested in the action of Lusztig's symmetries $T'_{1,2,\pm}$ on
the above weight spaces.

\lem{alge}
$T'_{1\pm}(F_1F_2\sv)=-v^{\pm(\mu_1+1)}F_1^{(\mu_1)}F_2\sv,\
T'_{1\pm}(F_2F_1\sv)=-v^{\pm\mu_1}F_1^{(\mu_1)}F_2\sv+F_2F_1^{(\mu_1)}\sv;\\
T'_{2\pm}(F_2F_1\sv)=-v^{\pm(\mu_2+1)}F_2^{(\mu_2)}F_1\sv,\
T'_{2\pm}(F_1F_2\sv)=-v^{\pm\mu_2}F_2^{(\mu_2)}F_1\sv+F_1F_2^{(\mu_2)}\sv;\\
T'_{2\pm}(F_2F_1^{(\mu_1)}\sv)=-v^{\pm(\mu_1+\mu_2)}F_2^{(\mu_1+\mu_2-1)}F_1^{(\mu_1)}\sv,\\
T'_{2\pm}(F_1^{(\mu_1)}F_2\sv)=-v^{\pm\mu_2}F_2^{(\mu_1+\mu_2-1)}F_1^{(\mu_1)}\sv+
F_1F_2^{(\mu_1+\mu_2-1)}F_1^{(\mu_1-1)}\sv;\\
T'_{1\pm}(F_1F_2^{(\mu_2)}\sv)=-v^{\pm(\mu_1+\mu_2)}F_1^{(\mu_1+\mu_2-1)}F_2^{(\mu_2)}\sv,\\
T'_{1\pm}(F_2^{(\mu_2)}F_1\sv)=-v^{\pm\mu_1}F_1^{(\mu_1+\mu_2-1)}F_2^{(\mu_2)}\sv+
F_2F_1^{(\mu_1+\mu_2-1)}F_2^{(\mu_2-1)}\sv;\\
T'_{1\pm}(F_1F_2^{(\mu_1+\mu_2-1)}F_1^{(\mu_1-1)}\sv)=
-v^{\pm(\mu_2+1)}F_1^{(\mu_2)}F_2^{(\mu_1+\mu_2-1)}F_1^{(\mu_1-1)}\sv,\\
T'_{1\pm}(F_2^{(\mu_1+\mu_2-1)}F_1^{(\mu_1)}\sv)=
-v^{\pm1}F_1^{(\mu_2)}F_2^{(\mu_1+\mu_2-1)}F_1^{(\mu_1-1)}\sv+
F_1^{(\mu_2-1)}F_2^{(\mu_1+\mu_2-1)}F_1^{(\mu_1)}\sv;\\
T'_{2\pm}(F_2F_1^{(\mu_1+\mu_2-1)}F_2^{(\mu_2-1)}\sv)=
-v^{\pm(\mu_1+1)}F_2^{(\mu_1)}F_1^{(\mu_1+\mu_2-1)}F_2^{(\mu_2-1)}\sv,\\
T'_{2\pm}(F_1^{(\mu_1+\mu_2-1)}F_2^{(\mu_2)}\sv)=
-v^{\pm1}F_2^{(\mu_1)}F_1^{(\mu_1+\mu_2-1)}F_2^{(\mu_2-1)}\sv+
F_2^{(\mu_1-1)}F_1^{(\mu_1+\mu_2-1)}F_2^{(\mu_2)}\sv.$
\elem

\prf We consider two subalgebras $\bU_1,\bU_2\subset\bU$ of type $A_1$: the
first one is generated by $E_1^{(r)},F_1^{(r)},\ r\in\BN$, the second one is
generated by $E_2^{(r)},F_2^{(r)},\ r\in\BN$. To prove the first formula, we
consider the $\bU_1$-submodule $M_1$ of $L(\Lambda)$ with the highest vector
$F_2\sv$ and canonical base $F_2\sv,F_1F_2\sv,\ldots,F_1^{(\mu_1)}F_2\sv,
F_1^{(\mu_1+1)}F_2\sv$. We also consider another $\bU_1$-submodule $M'_1$ of
$L(\Lambda)$ with the highest vector $\sw^+:=(v^{\mu_1}-v^{-\mu_1})F_1F_2\sv-
(v^{\mu_1+1}-v^{-\mu_1-1})F_2F_1\sv$ and the lowest vector (in the same canonical
base) $\sw^-=(v-v^{-1})F_1^{(\mu_1)}F_2\sv-(v^{\mu_1+1}-v^{-\mu_1-1})F_2F_1^{(\mu_1)}\sv$.
In effect, it is straightforward that $E_1\sw^+=0$, and it follows
from~\cite[Lemma~42.1.2.(d)]{lu} that $F_1\sw^-=0$; hence $F_1^{(\mu_1-1)}\sw^+=
a\sw^-$ for some $a\in\CA$. The fact that $a=1$ follows by comparing the
coefficients of $F_1F_2\sv$ in $\sw^+$ and of $F_1^{(\mu_1)}F_2\sv$ in $\sw^-$.
Now according to~\cite[Proposition~5.2.2.(a)]{lu}, $T'_{1\pm}\sw^+=\sw^-,\
T'_{1\pm}F_1F_2\sv=-v^{\pm(\mu_1+1)}F_1^{(\mu_1)}F_2\sv$. From this we deduce the
first two formulas. The other formulas are proved similarly. Say, to prove
the 5th and 6th formulas we consider the $\bU_2$-submodule $M_2$ of
$L(\Lambda)$ with the highest vector $F_1^{(\mu_1)}\sv$ and canonical base
$F_1^{(\mu_1)}\sv,F_2F_1^{(\mu_1)}\sv,\ldots,F_2^{(\mu_1+\mu_2-1)}F_1^{(\mu_1)}\sv,
F_2^{(\mu_1+\mu_2)}F_1^{(\mu_1)}\sv$. We also consider another $\bU_2$-submodule
$M'_2$ of $L(\Lambda)$ with the highest vector
$\sx^+:=(v^{\mu_1+\mu_2}-v^{-\mu_1-\mu_2})F_1^{(\mu_1)}F_2\sv-
(v^{\mu_2}-v^{-\mu_2})F_2F_1^{(\mu_1)}\sv$ and the lowest vector (in the same
canonical base)
$\sx^-=(v^{\mu_1+\mu_2}-v^{-\mu_1-\mu_2})F_1F_2^{(\mu_1+\mu_2-1)}F_1^{(\mu_1-1)}\sv-
(v^{\mu_1}-v^{-\mu_1})F_2^{(\mu_1+\mu_2-1)}F_1^{(\mu_1)}\sv$. Then $T'_{2\pm}\sx^+=\sx^-,\
T'_{2\pm}F_2F_1^{(\mu_1)}\sv=-v^{\pm(\mu_1+\mu_2)}F_2^{(\mu_1+\mu_2-1)}F_1^{(\mu_1)}\sv$.
From this we deduce the 5th and 6th formulas. And so on.
\epr

\ssec{topo}{Topology} Let $\BA_\BR$ be a 2-dimensional real vector space and let $\BA$ be its complexification with coordinates $(t_1,t_2)$
stratified by 3 lines: $t_1=0,\ t_2=0,\ t_1-t_2=0$. Let $\CL$ be the shriek
extension of the one-dimensional local system on the complement of the 3
lines with monodromies $v^{2\mu_1},v^{2\mu_2},v^{2\mu_3}$. In applications to algebra,
$2\mu_3=2$. The dual vector space $\BA^*$ has coordinates $(z_1,z_2)$, and
the dual stratification consists of the lines $z_1=0,\ z_2=0,\ z_1+z_2=0$.
This is the root hyperplane arrangement of type $A_2$.
There are 6 real chambers of this arrangement: $C_0$ is the dominant
chamber containing an interior point $z^{(e)}=(1,1)$; the other chambers
with interior points $z^{(1)}=(-1,2),\
z^{(21)}=(-2,1),\ z^{(121)}=z^{(212)}=(-1,-1),\ z^{(12)}=(1,-2),\ z^{(2)}=(2,-1)$.
The chambers are naturally numbered by the Weyl group $W$ of type $A_2$
generated by simple reflections $s_1,s_2$. For $w\in W$ we have
$C_w\ni z^{(w)}$, say $C_{s_1s_2s_1}$ ($C_{121}$ for short) contains $z^{(121)}$.
We have $6$ {\it real} affine lines $\ell_w,\ w\in W$, in $\BA_\BR$ given by equations
$z^{(w)}=1$. For example $\ell_e$, $\ell_1$ are given by the equations $t_1 + t_2 = 1$, $-t_1+2t_2 = 1$ respectively. More generally, for $\varepsilon\in \BR,\
\varepsilon>0$, let us denote by
$\ell_{w,\varepsilon}$ the real straight line given by the equation
$z^{(w)}=\varepsilon$.

The microlocalization (Fourier transform) $\bmu\CL$ is a certain
constructible complex on $\BA^*$. We will be interested only in its restriction to the complement of the 3 lines in $\BA^*$, which is a 2-dimensional local system. Let us describe this local system explicitely.

The stalk $\bmu^{(w)}\CL$ at $z^{(w)}$ equals the vanishing
cycles $\Phi_{z^{(w)}}\CL$. Let $i_w$ denote the inclusion
$\ell_w\hra \BA$. Then $\Phi_{z^{(w)}}\CL$ may be identified with
$H^1(\ell_w, i_w^*\CL)$.
It is a 2-dimensional vector space with the base
{\em dual}
to the basis $\ell'_w,\ell''_w$ of 1-cycles with coefficiens
in $i_w^*\CL^*$.

The 1-cycles are defined as follows:
$\ell'_e$ is the interval between the points $(1,0)$ and $(1/2,1/2)$;
$\ell''_e$ is the interval between the points $(1/2,1/2)$ and $(0,1)$;
$\ell'_1$ is the interval between the points $(1,1)$ and $(0,1/2)$;
$\ell''_1$ is the interval between the points $(0,1/2)$ and $(-1,0)$;
$\ell'_{21}$ is the interval between the points $(0,1)$ and $(-1/2,0)$;
$\ell''_{21}$ is the interval between the points $(-1/2,0)$ and $(-1,-1)$;
$\ell'_{121}$ is the interval between the points $(-1,0)$ and $(-1/2,-1/2)$;
$\ell''_{121}$ is the interval between the points $(-1/2,-1/2)$ and $(0,-1)$;
$\ell'_2$ is the interval between the points $(0,-1)$ and $(1/2,0)$;
$\ell''_2$ is the interval between the points $(1/2,0)$ and $(1,1)$;
$\ell'_{12}$ is the interval between the points $(-1,-1)$ and $(0,-1/2)$;
$\ell''_{12}$ is the interval between the points $(0,-1/2)$ and $(1,0)$;
$\ell'_{212}$ is the interval between the points $(-1,0)$ and $(-1/2,-1/2)$;
$\ell''_{212}$ is the interval between the points $(-1/2,-1/2)$ and $(0,-1)$.
Note that $\ell''_{212}=\ell''_{121},\ \ell'_{212}=\ell'_{121}$.
The dual basis to $\ell'_w,\ell''_w$ will be denoted by $\phi'_w,\phi''_w$
(in particular, $\phi''_{212}=\phi''_{121},\ \phi'_{212}=\phi'_{121}$).

More generally, for any $\varepsilon>0$ we have canonical isomorphisms
$$
\Phi_{z^{(w)}}\CL = H^1(\ell_w, i_w^*\CL) \iso
H^1(\ell_{w,\varepsilon}, i_{w,\varepsilon}^*\CL)
$$
where $i_{w,\varepsilon}: \ell_{w,\varepsilon}\hra \BA$, and we can define similar parallelly transported
bases in $H_1(\ell_{w,\varepsilon}, i_{w,\varepsilon}^*\CL^*)$.

For two neighbouring chambers $C_y,C_w,\ y,w\in W$, let $\gamma_{y,w}^\pm$ be a 
path from $C_y$ to $C_w$ obtained from a straight line interval modified near
the wall between these two chambers by going around it in 
the positive (resp. negative) imaginary
halfspace. We will keep the same notation for the induced operator (half
monodromy along $\gamma_{y,w}^\pm$) from $\Phi_{z^{(y)}}\CL$ to $\Phi_{z^{(w)}}\CL$.

\lem{topol}
$\gamma_{e,1}^\pm\phi'_e=-v^{\pm(\mu_1+\mu_3)}\phi''_1,\
\gamma_{e,1}^\pm\phi''_e=-v^{\pm\mu_1}\phi''_1+\phi'_1;\\
\gamma_{e,2}^\pm\phi''_e=-v^{\pm(\mu_2+\mu_3)}\phi'_2,\
\gamma_{e,2}^\pm\phi'_e=-v^{\pm\mu_2}\phi'_2+\phi''_2;\\
\gamma_{1,21}^\pm\phi'_1=-v^{\pm(\mu_2+\mu_1)}\phi''_{21},\
\gamma_{1,21}^\pm\phi''_1=-v^{\pm\mu_2}\phi''_{21}+\phi'_{21};\\
\gamma_{2,12}^\pm\phi''_2=-v^{\pm(\mu_1+\mu_2)}\phi'_{12},\
\gamma_{2,12}^\pm\phi'_2=-v^{\pm\mu_1}\phi'_{12}+\phi''_{12};\\
\gamma_{21,121}^\pm\phi'_{21}=-v^{\pm(\mu_3+\mu_2)}\phi''_{121},\
\gamma_{21,121}^\pm\phi''_{21}=-v^{\pm\mu_3}\phi''_{121}+\phi'_{121};\\
\gamma_{12,212}^\pm\phi''_{12}=-v^{\pm(\mu_3+\mu_1)}\phi'_{212},\
\gamma_{12,212}^\pm\phi'_{12}=-v^{\pm\mu_3}\phi'_{212}+\phi''_{212}.
$
\elem

\prf
All the formulas being similar, we prove the first two. For the transposed
map between dual spaces we must check that
$$
\gamma_{1,e}^\pm\ell'_1=\ell''_e,\
\gamma_{1,e}^\pm\ell''_1=-v^{\pm\mu_1}\ell''_e-v^{\pm(\mu_1+\mu_3)}\ell'_e.
$$
(Note that the second equality is equivalent to
$\gamma_{1,e}^\pm(\ell'_1+v^{\mp\mu_1}\ell''_1)=-v^{\pm\mu_3}\ell'_e$.)
To prove it, we rotate the line $\ell_1$
clockwise in $\BA_\BR$ with the point $(0,1/2)$ fixed and observe what
happens with the real cycles $\ell'_1$ and $\ell''_1$. At some critical moment
the rotated line becomes parallel to the $t_1$-axis, at this moment we must
pass for a short time into the complex upper (or lower) halfspace, and at the
end we get the line parallel to $\ell_e$. We see that
at the end of this rotation
$\ell'_1$ turns into $\ell''_e$, whereas $\ell''_1$ stretches and after
the critical moment turns into the necessary linear combination
of $\ell'_e$ and $\ell''_e$.
\epr

\rem{product}
Writing down the composition
$\gamma_{21,121}^\pm\circ\gamma_{1,21}^\pm\circ\gamma_{e,1}^\pm$ in our bases as
the product of matrices we find
$$\left(\begin{array}{cc}0&1\\
-v^{\pm(\mu_3+\mu_2)} &-v^{\pm\mu_3} \end{array}\right)\left(\begin{array}{cc}0&1\\
-v^{\pm(\mu_2+\mu_1)} &-v^{\pm\mu_2} \end{array}\right)\left(\begin{array}{cc}0&1\\
-v^{\pm(\mu_1+\mu_3)} &-v^{\pm\mu_1} \end{array}\right)=$$
$$\left(\begin{array}{cc}v^{\pm(\mu_1+\mu_2+\mu_3)}&0\\
0&v^{\pm(\mu_1+\mu_2+\mu_3)} \end{array}\right),$$
cf.~\cite[Corollary~5.9]{lu2}.
\erem

\rem{adjo}
In case $\mu_1=\mu_2=\mu_3=1$ all the six weight spaces considered
in~\refss{alg} coincide with $L(1,1)_{(0,0)}$ with the base
$F_1F_2\sv,\ F_2F_1\sv$. In this base the operator $T'_{1\pm}$ of the first line
of~\refl{alge} corresponding to the operator $\gamma_{e,1}^\pm$
of~\refl{topol} has the matrix $\left(\begin{array}{cc}-v^{\pm2}&-v^{\pm1}\\
0&1\end{array}\right)$, while the operator $T'_{2\pm}$ of the second line
of~\refl{alge} corresponding to the operator $\gamma_{e,2}^\pm$
of~\refl{topol} has the matrix $\left(\begin{array}{cc}1&0\\
-v^{\pm1}&-v^{\pm2}\end{array}\right)$. Note that
$(\gamma_{e,1}^-)^{-1}=\gamma_{e,1}^+$, and $(\gamma_{e,2}^-)^{-1}=\gamma_{e,2}^+$.
\erem

\ssec{disc}{Discussion}
We set $\mu_3=1$. The theory of factorizable sheaves~\cite{bfs} provides
a canonical isomorphism $\Phi_{z^{(e)}}\CL\simeq L(\Lambda)_{(\mu_1-1,\mu_2-1)}$.
The stalks of microlocalization at the other chambers $\Phi_{z^{(w)}}\CL$
{\em do not} have an algebraic interpretation in the framework of this
theory.\footnote{see ~\refco{refor} however.}
However, the comparison of~\refl{alge} and~\refl{topol} shows that
the monodromy of the local system $\bmu\CL$ (as the automorphism group of
$\Phi_{z^{(e)}}\CL\simeq L(\Lambda)_{(\mu_1-1,\mu_2-1)}$) can be expressed in terms
of Lusztig's symmetries $T'_{1,2\pm},T''_{1,2\pm}$. In fact, the comparison
of~\refl{alge} and~\refl{topol} suggests a much more precise relation,
in particular, between a natural topological basis in $\Phi_{z^{(e)}}\CL$ and the
canonical basis on the algebraic side. Unfortunately, we have no clue how to
define such a topological basis in general. However, the relation between
the monodromy and Lusztig's symmetries seems to generalize. This is the subject
of the main body of the note.

\sec{Nota}{Coxeter categories}

\ssec{quadr}{Notations}
Let us set up a few notations related to a simple Lie algebra $\fg$ with Cartan
subalgebra $\fh$ and Borel subalgebra $\fh\subset\fb\subset\fg$.
The set of simple coroots is denoted by $I$; for $i\in I$ the corresponding
simple coroot is denoted $\check\alpha_i$ or sometimes simply $i$. The
corresponding simple root is denoted $\alpha_i$ or sometimes $i'$.
We fix a Weyl group invariant symmetric bilinear form $?\cdot?$ on $\fh^*$
such that the square length of a {\em short} root is $\alpha_i\cdot\alpha_i=2$.
This bilinear form gives rise to an isomorphism $\fh\iso\fh^*$ so that
the coroot lattice $Y$ generated by $\{\check\alpha_i\}_{i\in I}$ embeds into
$\fh^*$. We then have $\check\alpha_i\cdot\check\alpha_i\in\{2,1,\frac{2}{3}\}$,
and $\alpha_i\cdot\alpha_i\in\{2,4,6\}$.
We set $d_i=\alpha_i\cdot\alpha_i/2$. Let $d$ be the ratio of the square
lengths of the long and short roots, so that $d\in\{1,2,3\}$. We set
$\check{d}_i=d/d_i$. Then
$$
\langle\alpha_i,\check\alpha_j\rangle=
\frac{\alpha_i\cdot\alpha_j}{d_j}=d_i\check\alpha_i\cdot\check\alpha_j=
d\frac{\check\alpha_i\cdot\check\alpha_j}{\check{d}_i}.
$$

\ssec{sal}{The fundamental groupoid of $\fh_{D'}^\reg$}
We follow the notations of~\cite{at}. Let $D$ be the Dynkin diagram
of the simple Lie algebra $\fg$ with Cartan $\fh$ (so that $I$ is the
set of vertices of $D$). The root system of
$\fh\subset\fg$ is $R_D\subset\fh^*$. The complement in $\fh$ to the root
hyperplanes is the open subset $\fh^{\on{reg}}$.

Given a subset $D'$ of the set
of vertices of $D$, we denote by $\fh_{D'}$ the quotient of $\fh$
by the center of the corresponding Levi subalgebra $\fl_{D'}\subset\fg$. In other words $\fh_{D'}^*\subset \fh_{D}^*$
is spanned by the simple roots corresponding to the vertices from
$D'$.
We denote by $\fh_{D'}^{\on{reg}}$ the complement
in $\fh_{D'}$ to the root hyperplanes of the root subsystem $R_{D'}$ corresponding to $D'\subset D$.

We recall the Salvetti presentation of the fundamental groupoid of
$\fh_{D'}^\reg$, cf.~\cite{sa}. Let $\fh_{D',\BR}^\reg$ denote the set of real points of
$\fh_{D'}^\reg$. It is a union of the connected components called
{\em chambers}. We fix a chamber $C^{D'}_0$ formed by the points with
{\em positive} coordinates (in the basis of fundamental coweights). The Weyl
group $W_{D'}$ acts on the set $\bC^{D'}$ of chambers simply transitively
on the left. The choice of $C^{D'}_0$ identifies $\bC^{D'}$ with $W_{D'}$, and
defines the {\em right} action of $W_{D'}$ on $\bC^{D'}$ (transferred from the
right action of $W_{D'}$ on itself). The set of walls of $C^{D'}_0$ is
canonically identified with the set of vertices of $D'$. The left action of
$W_{D'}$ on $\bC^{D'}$ extends this identification to any chamber. For $i\in D'$
the right action of a simple reflection works as follows: $C\cdot s_i$ is a
unique neighbouring chamber $C'$ having the $s_i$-wall in common with $C$.

The set of objects of the fundamental groupoid $\Pi(\fh_{D'}^\reg)$ is $\bC^{D'}$.
Given a straight line interval $\gamma$ connecting the endpoints
$\gamma_1\in C_1$ and $\gamma_2\in C_2$ and intersecting only one wall at a
time, we define the morphisms $\gamma^\pm\in\on{Mor}_{\Pi(\fh_{D'}^\reg)}(C_1,C_2)$
as follows.
The path $\gamma^+$ (resp. $\gamma^-$) coincides with $\gamma$ away from the
small neighbourhoods of its intersection with walls, where $\gamma^+$
(resp. $\gamma^-$) goes around the intersection in the positive
(resp. negative) imaginary direction in $\fh_{D'}^\reg$. According to Salvetti,
$\Pi(\fh_{D'}^\reg)$ is generated by the set of morphisms $\gamma^\pm$ with
relations $\beta^\pm=\gamma^\pm$ provided $\gamma_1,\beta_1$ lie in the same
chamber $C_1$, and $\gamma_2,\beta_2$ lie in the same chamber $C_2$.

\ssec{spec}{The fundamental groupoid of $N_{\fh_{D'/D''}/\fh_{D'}}^\reg$}
Given a third subdiagram $D''\subset D'\subset D$, we consider the linear subspace
$\fh_{D'/D''}\subset\fh_{D'}$ spanned by the {\em fundamental coweights} in
$D'-D''\subset D'$. For example, $\fh_{D/D''}\subset\fh_D$ is the center of the Levi $\fl_{D''}\subset\fg$. We have an exact sequence
$$
0\lra \fh_{D'/D''} \lra \fh_{D'} \lra \fh_{D''}\lra 0
$$
which may serve as another definition of $\fh_{D'/D''}$.

We denote by $\fh_{D'/D''}^\reg$ the complement in $\fh_{D'/D''}$
to the root hyperplanes (roots in $R_{D'}$) not containing $\fh_{D'/D''}$.
The connected components of the real part $\fh_{D'/D'',\BR}^\reg$ are called
{\em chambers}; the set of chambers is denoted $\bC^{D'/D''}$. It is
naturally isomorphic to the set of parabolics in $\fg_{D'}$ containing the
standard Levi $\fl_{D''}$, see e.g.~\cite[I.1.10]{mw}.
We say that a
chamber $C\in\bC^{D'}$ is {\em adjacent} to $\fh_{D'/D''}\subset\fh_{D'}$
if the intersection of the closure $\overline{C}$ with $\fh_{D'/D''}$ has
the maximal (real) dimension $\dim\fh_{D'/D''}$; then this intersection is
the closure of a chamber in $\bC^{D'/D''}$ to be denoted $\pi(C)$. The set
of chambers adjacent to $\fh_{D'/D''}$ is denoted $\bA^{D'/D''}$. Thus
we have a projection $\pi\colon \bA^{D'/D''}\to\bC^{D'/D''}$.

The natural projection $\on{pr}\colon \fh_{D'}\to\fh_{D''}$ (see~\refss{sal})
works in the bases
of fundamental coweights as follows: $\on{pr}\check\omega_i=0$ if
$i\in D'-D''\subset D'$; and if $i\in D''\subset D'$, then $\check\omega_i$
goes to the corresponding fundamental coweight in $\fh_{D''}$. Given a 
chamber $C\in\bA^{D'/D''}$, its projection $\on{pr}C$ is a chamber in
$\bC^{D''}$. Thus we have a projection $\on{pr}\colon \bA^{D'/D''}\to\bC^{D''}$.

\lem{adj}
The product $\on{pr}\times\pi\colon \bA^{D'/D''}\to\bC^{D''}\times\bC^{D'/D''}$
establishes a one-to-one correspondence. $\square$
\elem

\defe{pipi}
(a) For a chamber $C\in\bC^{D''}$, we define a subgroupoid
$\Pi_{\on{pr}^{-1}(C)}(\fh_{D'}^\reg)\subset\Pi(\fh_{D'}^\reg)$ as follows: the objects are
$\on{pr}^{-1}(C)\subset\bA^{D'/D''}\subset\bC^{D'}$, and the morphisms are
generated by $\gamma^\pm$ where $\gamma$ is a straight line interval
{\em parallel} to $\fh_{D'/D''}$ (i.e.\ such that $\on{pr}\gamma$ is a point).

(b) For a chamber $C\in\bC^{D'/D''}$, we define a subgroupoid
$\Pi_{\pi^{-1}(C)}(\fh_{D'}^\reg)\subset\Pi(\fh_{D'}^\reg)$ as follows: the objects are
$\pi^{-1}(C)\subset\bA^{D'/D''}\subset\bC^{D'}$, and the morphisms are generated
by $\delta^\pm$ where $\delta$ is a straight line interval inside the union
of closures of chambers adjacent to $C$.

\begin{figure}
$$\xymatrix{
\ar@{-}[ddddddddrrrrrrrr] & & & & 
\ar@{--}[dddddddd]^(.00){\fh_{D''}}_(.10){C_N}_(.90){C_S} & & & & \\
& & & & & & & & \\
& C_{NW}\ar@{.}[rrrrrr]_(.60){\gamma_N} \ar@{.}[dddd]^(.60){\delta_W} 
& & & & & & C_{NE}\ar@{.}[dddd]_(.60){\delta_E} & \\
& & & & & & & & \\
\ar@{-}[rrrrrrrr]_(1.00){\fh_{D'/D''}}^(.93){C_E}^(.07){C_W} & & & & & & & & \\
& & & & & & & & \\
& C_{SW}\ar@{.}[rrrrrr]^(.60){\gamma_S} & & & & & & C_{SE} & \\
& & & & & & & & \\
\ar@{-}[uuuuuuuurrrrrrrr] & & & & & & & & 
}$$
\caption{Topographical example}
\label{topograph}
\end{figure}
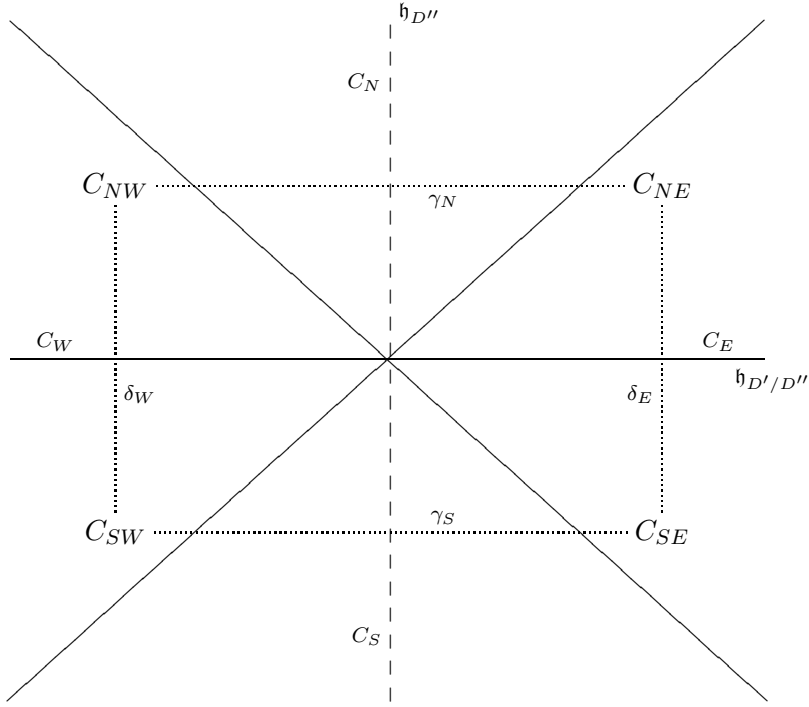

(c) A subgroupoid $\Pi_{\bA^{D'/D''}}(\fh_{D'}^\reg)\subset\Pi(\fh_{D'}^\reg)$
is generated by all the groupoids in (a,b) above. That is, its objects are
$\bA^{D'/D''}$, and the morphisms are all the possible products of morphisms
in (a,b) above (see Figure~\ref{topograph}).
\edefe

\lem{pip}
(a) For any $C\in\bC^{D''},\ \pi$ induces an equivalence
$\Pi_{\on{pr}^{-1}(C)}(\fh_{D'}^\reg)\iso\Pi(\fh_{D'/D''}^\reg)$.

(b) For any $C\in\bC^{D'/D''},\ \on{pr}$ induces an equivalence
$\Pi_{\pi^{-1}(C)}(\fh_{D'}^\reg)\iso\Pi(\fh_{D''}^\reg)$.

(c) 
The natural projection morphisms 
$\pi\colon \Pi_{\bA^{D'/D''}}(\fh_{D'}^\reg)\to\Pi(\fh_{D'/D''}^\reg$ and
$\on{pr}\colon \Pi_{\bA^{D'/D''}}(\fh_{D'}^\reg)\to\Pi(\fh_{D''}^\reg)$ give rise
to an equivalence $\Pi_{\bA^{D'/D''}}(\fh_{D'}^\reg)\iso
\Pi(\fh_{D'/D''}^\reg)\times\Pi(\fh_{D''}^\reg)$.
\elem

\prf (M.~Kapranov)
The relations in the Salvetti complex~\cite{sa} follow from a cell
decomposition of the complement which is glued out of intervals,
$2n$-gons (for any codim 2 cell where $n$ hyperplanes meet) and so on,
and the relations in the fundamental groupoid are obtained from the
2-skeleton i.e., from these $2n$-gons.

So the 2-dimensional case implies the general one. See example with 
$n=3$ in Figure~\ref{topograph}.
We keep ``vertical" and ``horizontal'' $2k$-gons (contributing to
$\Pi_{\pi^{-1}(C)}(\fh_{D'}^\reg)$ or $\Pi_{\on{pr}^{-1}(C)}(\fh_{D'}^\reg)$)
intact, and replace the remaining $2n$-gons with rectangles like the
dotted one in Figure~\ref{topograph}. It follows that the ``horizontal''
and ``vertical'' morphisms commute.
This produces a 2-dimensional CW-subcomplex of the complement
which is (the 2-skeleton of) the product of two separate 2-dimensional
subcomplexes.
\epr

\ssec{speci}{Specialization} We consider the normal bundle
$N_{\fh_{D'/D''}/\fh_{D'}}\iso \fh_{D'/D''}\times \fh_{D''}$, and
its open subspace
$N_{\fh_{D'/D''}/\fh_{D'}}^\reg=\fh_{D'/D''}^\reg\times\fh_{D''}^\reg
\subset  N_{\fh_{D'/D''}/\fh_{D'}}$ with Poincar\'e groupoid
$\Pi_{\bA^{D'/D''}}(\fh_{D'}^\reg)$. Then the Verdier specialization of a local system
on $\fh_{D'}^\reg$ along $\fh_{D'/D''}$ will be a well defined local system
on $N_{\fh_{D'/D''}/\fh_{D'}}^\reg$. At the level of representations of Poincar\'e
groupoids, the specialization is nothing but restriction to
$\Pi_{\bA^{D'/D''}}(\fh_{D'}^\reg)$.

\defe{Coxeter}
(A) A {\em pure Coxeter category of type $D$} is the collection of the 
following data:

(a) A category $\CalC_{D'}$ for any subset $D'\subset D$;

(b) For $D''\subset D'$ a local system of restriction functors
$F_{D'D''}\colon \CalC_{D'}\to\CalC_{D''}$ on $\fh_{D'/D''}^\reg$;

(c) For $D'''\subset D''\subset D'$ an isomorphism of local systems of functors:
$$
\phi_{D'D''D'''}\colon \on{Sp}_{\fh_{D'/D''}}F_{D'D'''}\iso F_{D''D'''}\circ F_{D'D''}
$$
on $\fh_{D'/D''}^\reg\times\fh_{D''/D'''}^\reg$ which satisfy the natural ``cocycle''
or ``pentagon" identity associated with
$D^{iv}\subset D'''\subset D''\subset D'$.

(d) In case $D'''$ is disjoint from $D'$ (i.e. no vertex of $D'''$ is connected
by an edge to a vertex of $D'$, and $D'''\cap D'=\emptyset$) we have a 
canonical isomorphism $\fh_{D'/D''}^\reg=\fh_{(D'\cup D''')/(D''\cup D''')}^\reg$,
and we are given a homomorphism of local systems of endomorphism algebras
$\eta\colon \End(F_{D'D''})\to\End(F_{(D'\cup D''')(D''\cup D''')})$.

(B) A {\em tensor Coxeter category of type $D$} is the additional datum of
braided balanced tensor structures on $\CalC_{D'}$ such that

(i) The pullback of $F_{D'D''}$ to the universal cover 
$\widetilde\fh{}^\reg_{D'/D''}$ is a (trivial) local system of functors
equipped with tensor structures
$\widetilde{F}_{D'D''}\colon \CalC_{D'}\to\CalC_{D''}$.

(But we {\em do not} require them to respect the balance and braiding. 
Also, the monodromy isomorphisms of stalks
$\gamma_*\colon (F_{D'D''})_x \iso (F_{D'D''})_y$ where
$\gamma\colon x\lra y$ is a path in $\Pi(\fh_{D'/D''}^\reg)$
{\em are not} required to be morphisms of tensor functors.
Neither are the isomorphisms (Ad) above required to respect the tensor
structure.);

(ii) The isomorphisms of (c) above pulled back to 
$\widetilde\fh{}^\reg_{D'/D''}\times\widetilde\fh{}^\reg_{D''/D'''}$ are
isomorphisms of tensor functors;

(iii) Let $\gamma_0\in \pi_1(\fh_{D'/D''}^\reg)$ be the generator
of the centre; geometrically it is a loop $\exp(2\pi i\theta)\cdot x,\
0\leq\theta\leq 1,\ x\in\fh^\reg_{D'/D''}$. The automorphism
$\gamma_{0*}\colon F_{D'D''}\iso F_{D'D''}$ (it is the automorphism
induced by the $\BC^*$-monodromic structure on the sheaf $F_{D'D''}$)
is inverse
to the ratio of the balance automorphisms of the identity functors of
$\CalC_{D''}$ and $\CalC_{D'}$.\footnote{We thank A.~Appel and
V.~Toledano Laredo for 
correcting mistakes in the original version of the definition.}

\edefe

\ssec{comp}{Comparison with the Appel-Toledano-Laredo Coxeter braided 
tensor categories}
If we impose an additional assumption that the local systems 
in~\refd{Coxeter}(Ab)
are lifted from $\fh_{D'/D''}^\reg/W_{D'/D''}$ (quotient with respect to the free
action of the finite group $W_{D'/D''}:=\on{Norm}_{L_{D'}}(L_{D''})/L_{D''}$
(normalizer of the Levi
subgroup $L_{D''}$ in $L_{D'}\subset G$, modulo $L_{D''}$)),
then we get an equivalent version of~\cite[Definitions~3.10,~4.1]{at}.

In the example of factorizable sheaves $\CF\CS_D$ (\refss{Coxeterfs} below),
the balance on an irreducible sheaf $\CL(\lambda)$ is multiplication
by $\zeta^{\lambda\cdot(\lambda+2\rho)}$. Factorizable sheaves $\CF\CS_\emptyset$
for Levi=Cartan also have a nontrivial braiding and balance; namely, on an
irreducible sheaf $\CL_\emptyset(\mu)$ the balance is multiplication by
$\zeta^{\mu\cdot(\mu+2\rho)}$. The ratio of these two balances on a
weight component $\CL(\lambda)^\alpha$ of $\CL(\lambda)$ is
$\zeta^{\lambda\cdot(\lambda+2\rho)-(\lambda-\alpha)\cdot(\lambda-\alpha+2\rho)}$ and
coincides with the monodromy automorphism of the monodromic sheaf
$\CL(\lambda)^\alpha$.

The identity $\Delta_i(T_i)=R_i^{21}\cdot(T_i\otimes T_i)$, and
more generally, for any $D'\subset D,\
\Delta_{D'}(T_{w_0^{D'}})=R_{D'}^{21}\cdot(T_{w_0^{D'}}\otimes T_{w_0^{D'}})$, implies
$\Delta_{D'}(T_{w_0^{D'}})^2(T^2_{w_0^{D'}}\otimes T^2_{w_0^{D'}})^{-1}=
R_{D'}^{12}\circ R_{D'}^{21}$,
which in view of~\refd{Coxeter}(Biii) is nothing but the usual relation between
the braiding and the balance.

\sec{algeb}{Algebra}

\ssec{ti}{The Lusztig symmetries}
Given $\zeta\in\BC,\ \zeta^6\ne1$, we consider the Lusztig small quantum group $\bu_{D'}$ (see e.g.~\cite[0.2.12]{bfs}). We extend it by the projectors to the weight spaces $1_\lambda,\ \lambda\in X$, to obtain the algebra $\dot\bu_{D'}$ such that $\on{Rep}(\dot\bu_{D'})=\CalC_{D'}$ (notations of~\cite[0.2.11-0.2.13]{bfs}). The algebra $\dot\bu_{D'}$ is a subalgebra of the Lusztig big quantum group
$_R\dot\bU_{D'}$~\cite[Chapter~31]{lu} (where $R\colon \BZ[v^{\pm1}]\to\BC,\ v\mapsto\zeta$), generated by $E_i=E_i^{(1)},\ F_i=F_i^{(1)},\ i\in D'$, and $1_\lambda,\ \lambda\in X$. According to~\cite[Chapters~33,35]{lu}, if $\zeta$ is a root
of unity (primitive of order $\ell$), there is a reductive
algebraic group $\check{G}_{D',\zeta}$ with Cartan torus 
$\check{T}_{D',\zeta}\subset\check{G}_{D',\zeta}$ and a tensor functor 
$\on{Fr}^*\colon \on{Rep}(\check{G}_{D',\zeta})\to\on{Rep}(_R\dot\bU_{D'})$
(pullback with respect to the quantum Frobenius homomorphism).
Note that the character lattice $X^*(\check{T}_{D',\zeta})$ is naturally a
sublattice of the weight lattice $X$.

The Lusztig symmetries $T'_{i,e},T''_{i,e},\ i\in D',\ e=\pm1$, of $_R\dot\bU_{D'}$~\cite[41.1.8]{lu} clearly preserve the subalgebra
$\dot\bu_{D'}$ and restrict to the same named symmetries of this subalgebra.
We define a functor $T_\bu$ from $\Pi(\fh_{D'}^\reg)$ to the category of $\BC$-algebras on generators as follows: $T_\bu(C)=\dot\bu_{D'}$ for any $C\in\bC$; for $\gamma$ a straight line interval connecting the endpoints in two neighbouring chambers $C_1,C_2$ with the common wall of type $s_i,\ i\in D'$, we set $T_\bu(\gamma^+)=T'_{i,1}$ (resp. $T''_{i,-1}$) and  $T_\bu(\gamma^-)=T'_{i,-1}$ (resp. $T''_{i,1}$), if $\gamma$ goes from a
 Bruhat smaller chamber to the bigger one (resp. from a Bruhat bigger chamber
 to the smaller one). According to~\cite[Theorem~39.4.3]{lu}, $T_\bu$ is well defined.

Given an integrable $_R\dot\bU_{D'}$-module $M$ with the Lusztig symmetries $T'_{i,e},\ T''_{i,e}\colon M\to M$~\cite[41.2.3]{lu}, we define a functor $T_M$ from $\Pi(\fh_{D'}^\reg)$ to the category of $\BC$-vector spaces on generators as follows: $T_M(C)=M$ for any $C\in\bC$; for $\gamma$ a straight line interval connecting the endpoints in two neighbouring chambers $C_1,C_2$ with the common wall of type $s_i,\ i\in D'$, we set $T_M(\gamma^+)=T'_{i,1}$ (resp. $T''_{i,-1}$) and  $T_M(\gamma^-)=T'_{i,-1}$ (resp. $T''_{i,1}$), if $\gamma$ goes from a Bruhat smaller chamber to the bigger one (resp. from a Bruhat bigger chamber to the smaller one). According to~\cite[Proposition~41.2.4]{lu}, $T_M$ is well defined.

Let $_R\CalC_{D'}$ (resp. $\CalC_{D'}$) denote the category of integrable
$_R\dot\bU_{D'}$-modules (resp. $\dot\bu_{D'}$-modules),
and let $\Upsilon\colon _R\CalC_{D'}\to\CalC_{D'}$ stand for the restriction functor. In the previous paragraph we have defined the local system $_RF_{D',\emptyset}$
on $\fh_{D'}^\reg$
of restriction functors $_R\CalC_{D'}\to\on{Vect}_X=\ _R\CalC_{\emptyset}$ to the category of $X$-graded $\BC$-vector spaces.

\prop{ups}
There exists a unique local system $F^\CalC_{D',\emptyset}$ on $\fh_{D'}^\reg$ of
restriction functors $\CalC_{D'}\to\on{Vect}_X=\CalC_{\emptyset}$ such that
$_RF_{D',\emptyset}=F^\CalC_{D',\emptyset}\circ\Upsilon$.
\eprop

\prf
According to~\cite[Theorem~4.7]{ag}, we view $\CalC_{D'}$ as the category
of Hecke-eigen-objects in $_R\CalC_{D'}$. That is, an object of $\CalC_{D'}$ is
an object $M$ of $_R\CalC_{D'}$ endowed with a collection of isomorphisms
$\alpha_V\colon \on{Fr}^*(V)\otimes M\iso
\on{Res}^{\check{T}_{D',\zeta}}_{\check{G}_{D',\zeta}}(V)\otimes M,\ 
V\in\on{Rep}(\check{G}_{D',\zeta})$. 
Since the Lusztig symmetries act on $\on{Fr}^*(V)$ and on 
$\on{Res}^{\check{T}_{D',\zeta}}_{\check{G}_{D',\zeta}}(V)$, a Hecke-eigen-object 
$(M,\alpha)$ gives rise to a representation $T_{(M,\alpha)}$ of 
$\Pi(\fh^\reg_{D'})$. Hence the action of $\Pi(\fh_{D'}^\reg)$ on
$_RF_{D',\emptyset}$ canonically extends to the action of $\Pi(\fh_{D'}^\reg)$ on
$F_{D',\emptyset}$.
\epr

\rem{restr}
For example, if $\lambda\in X$ is a dominant $\ell$-restricted weight
(recall that $\ell$ is the order of $\zeta$), then the irreducible 
$\dot\bu_{D'}$-module $L^\lambda_\bu$ with highest weight $\lambda$ is the 
restriction of the irreducible $_R\dot\bU_{D'}$-module $L^\lambda_\bU$
with highest weight $\lambda$, and $T_{L^\lambda_\bu}=T_{L^\lambda_\bU}$.
\erem

\ssec{Coxeterc}{A Coxeter structure on $_R\CalC,\ \CalC$}
We need to define the local systems of restriction functors $F^\CalC_{D'D''}$,
not just $F^\CalC_{D',\emptyset}$ as in the previous subsection.
To this end we restrict the
action of the fundamental groupoid $\Pi(\fh_{D'}^\reg)$ defined in~\refp{ups}
to the subgroupoid
$\Pi_{\on{pr}^{-1}(C_0^{D''})}(\fh_{D'}^\reg)\subset\Pi(\fh_{D'}^\reg)$,
see~\refd{pipi}(a). More precisely, we consider a wall between two neighbouring chambers
$c_1,c_2$ of $\fh_{D'/D'',\BR}^\reg$. Let $C_1$ (resp. $C_2$) be a (unique) chamber of $\fh_{D',\BR}^\reg$ adjacent to $c_1$ (resp. $c_2$) such that
$\on{pr}(C_1)=\on{pr}(C_2)=C_0^{D''}$ (the fundamental chamber of $\fh_{D'',\BR}^\reg$). Let $\gamma$ be a straight line interval going from
$c_1$ to $c_2$, and let $\Gamma$ be its lift going from $C_1$ to $C_2$
parallelly to $\fh_{D'/D''}$. In the notations of~\refss{ti} we set
$T_M^{D'/D''}(\gamma^+)=\on{Id}_M,\
T_M^{D'/D''}(\gamma^-)=T_M^{-1}(\Gamma^-)\circ T_M(\Gamma^+)$.

\lem{commut}
The action of $\Pi(\fh_{D'/D''}^\reg)$
on a $_R\dot\bU_{D'}$-module $M$ commutes with the action of
$_R\dot\bU_{D''}\subset\ _R\dot\bU_{D'}$.
\elem

\prf
Note that (since $c_1$ and $c_2$ are neighbours) the types of walls intersected
by $\Gamma$ are all in $D''$ except for exactly one $d\in D'-D''$. Let $C'_1$
be the chamber adjacent to $c_1$ such that $\on{pr}(C'_1)=w_0^{D''}C_0^{D''}$.
Then $C_2=w_0^{D''\sqcup d}C'_1$. Let $\Delta$ be a straight line interval
going from $C'_1$ to $C_1$, and ending at the starting point of $\Gamma$, and
let $\Gamma\Delta$ be the concatenation of $\Gamma$ and $\Delta$.
Then $T_M^{-1}(\Gamma^-)\circ T_M(\Gamma^+)=T_M^{-1}(\Delta^+)\circ T_M(\Delta^-)\circ T_M^{-1}(\Gamma\Delta^-)\circ T_M(\Gamma\Delta^+)$.
It suffices to prove that $T_\bU^{-1}(\Gamma^-)\circ T_\bU(\Gamma^+)=
\on{Id}_{_R\dot\bU_{D''}}$. According to~\cite[Corollary~5.9]{lu2} or~\cite{kt},
we have
$T_\bU^{-1}(\Gamma\Delta^-)\circ T_\bU(\Gamma\Delta^+)(E_i)=\tilde{K}{}_i^{-2}E_i,\
T_\bU^{-1}(\Gamma\Delta^-)\circ T_\bU(\Gamma\Delta^+)(F_i)=F_i\tilde{K}{}_i^2$ for any
$i\in D''\sqcup d$, and
$T_\bU^{-1}(\Delta^-)\circ T_\bU(\Delta^+)(E_i)=\tilde{K}{}_i^{-2}E_i,\
T_\bU^{-1}(\Delta^-)\circ T_\bU(\Delta^+)(F_i)=F_i\tilde{K}{}_i^2$ for any $i\in D''$.
\epr

Now by~\refl{pip}(a) (and~\refp{ups}) we obtain the desired local system
of restriction functors $_RF_{D'D''}\colon _R\CalC_{D'}\to\ _R\CalC_{D''}$
(resp. $F^\CalC_{D'D''}\colon \CalC_{D'}\to\CalC_{D''}$) on $\fh_{D'/D''}^\reg$.
The isomorphisms of~\refl{pip}(c) give rise to the isomorphisms of~\refd{Coxeter}(Ac).
The conditions of~\refd{Coxeter}(B) are satisfied trivially.

\rem{tole}
The Coxeter structure on $_R\CalC$ studied in~\cite{tr} differs from ours by
the twist by an invertible local system. More precisely, for a weight
component $M_\lambda\subset M$, in the setup of~\refss{ti}, the Coxeter
structure of~\cite[4.1.3]{tr}
$T_{M_\lambda}^{TL}(\gamma^+)=\zeta^{d_i\langle\check\alpha_i,\lambda\rangle^2/4}T''_{i,+1}=
(-1)^{\langle\check\alpha_i,\lambda\rangle}\zeta^{d_i\langle\check\alpha_i,\lambda\rangle}
\zeta^{d_i\langle\check\alpha_i,\lambda\rangle^2/4}T'_{i,+1}=
(-1)^{\langle\check\alpha_i,\lambda\rangle}
\zeta^{d_i\langle\check\alpha_i,\lambda\rangle+d_i\langle\check\alpha_i,\lambda\rangle^2/4}
T_{M_\lambda}(\gamma^+)$
(the second equality is~\cite[5.2.3.b)]{lu}) for $\gamma$ going through an $s_i$-wall from a Bruhat smaller chamber to a Bruhat bigger one;
$T_{M_\lambda}^{TL}(\gamma^-)=\zeta^{d_i\langle\check\alpha_i,\lambda\rangle^2/4}T''_{i,+1}=
\zeta^{d_i\langle\check\alpha_i,\lambda\rangle^2/4}T_{M_\lambda}(\gamma^-)$
for $\gamma$ going through an $s_i$-wall from a Bruhat bigger chamber to a Bruhat smaller one; the remaining two halfmonodromies are the inverses of the above two.

Note that if $s_i\lambda=\lambda$ then the scalar factors above are identically
equal to one. We define $M_{W\lambda}:=\bigoplus_{\mu\in W\lambda}M_\mu$, the direct
sum over the Weyl group orbit of $\lambda$. Since $T_{M_{W\lambda}}^{TL}(\gamma^\pm)$
arise from a local system on $\fh^\reg/W$~\cite{tr} (i.e. $_RF_{D\emptyset}^{TL}$
possesses a $W$-equivariant structure), it follows that $_RF_{D\emptyset}$ also
possesses a $W$-equivariant structure.
\erem

\sec{topolo}{Topology}

\ssec{err}{Erratum to~\cite{bfs}}
We take this opportunity to correct a blunder pertaining to
the non-simply laced case of ~\cite{bfs}.
Let us define the quantities $i'\cdot j'$ as $\alpha_i\cdot\alpha_j$
in the sense of~\refss{quadr}.
Then throughout~\cite{bfs} in all formulas
the occurences of $i\cdot j$ should be replaced by $i'\cdot j'$.

For example:

--- in~\cite[Part 0, 2.1]{bfs}: $\langle i, j'\rangle = 2i'\cdot j'/i'\cdot i'$, and
$d_i = i'\cdot i'/2$;

--- in the relation~\cite[Part 0, 2.7(d)]{bfs} one should
replace $\zeta^{i\cdot j}$ by $\zeta^{i'\cdot j'}$,

Thus, in~\cite[Part 0, 2.7]{bfs} $\tilde{K}_i=K_i^{d_i}$ as before,
and $\zeta_i=\zeta^{d_i}$,
but the meaning of $d_i$ should be changed:
 $d_i$ is defined
{\em not} as half the square length of the coroot $\check\alpha_i$ but as half the square length of the root $\alpha_i$.

--- On the geometric side, the monodromy of the cohesive local
system corresponding to a full counterclockwise turn of a point $i$ around $j$
should be $\zeta^{-2i'\cdot j'}$, cf. ~\cite[Part 0, 3.10]{bfs}, and~\refss{rev1}
below.

To summarize, the main assertion of ~\cite{bfs}
(reviewed below in more details) consists of two
parts: first, an equivalence of the geometric category $\CF\CS$
with a category of graded modules over the algebra $\fu$ defined
in~\cite[Part 0, 2.7]{bfs}. This assertion is true, and our correction just
replaces the root system by the dual one on both sides.
The second assertion is an identification of $\fu$ with the
Lusztig's small quantum group. This identification is described
in~\cite[Part 0, 2.12, or Part II, 12.5]{bfs} and should be corrected:
the ``geometric" algebra $\fu$ is isomorphic to the ``Langlands dual"
Lusztig's algebra connected with the {\em dual} root system.

\vspace{.2cm}

This replacement of the root system by its dual is a rather subtle point. Its
origin lies in the definition of the braiding in~\cite{lu}, cf.
the proof of~\cite[Lemma~32.2.3]{lu}.

\vspace{.2cm}

Also, there is a misprint in the definition of a balance
in~\cite[IV.6.6]{bfs}: $n(\lambda)$ must be replaced by
$2n(\lambda)=\lambda\cdot(\lambda+2\rho)$.

\ssec{rev1}{A review of~\cite{bfs}: cohesive system and algebra $\fu^-$}
For $\beta\in\BN[I]$ we consider the configuration space $\BA^\beta$ of colored
divisors on the complex affine line $\BA^1$. The open subspace
$\overset{\circ}\BA{}^\beta\subset\BA^\beta$ of multiplicity free divisors carries
a 1-dimensional {\em cohesive local system} $\CJ^\beta$
with the following monodromies:
$\zeta^{-2\alpha_i\cdot\alpha_j}$ when a point of colour $i$ goes
counterclockwise around a point of colour $j\ne i$, and
$-\zeta^{-\alpha_i\cdot\alpha_i}$ when two points of colour $i$ trade their
positions going around a halfcircle counterclockwise. We denote by
$\CI^\beta$ the Goresky-MacPherson extension of $\CJ^\beta$ to $\BA^\beta$
(a perverse sheaf). Given two disjoint open discs
$\BA^1\supset D(p_i,\varepsilon_i),\ i=1,2$,
with centers in $p_i$ of radii $\varepsilon_i$, and $\beta_{1,2}\in\BN[I]$,
we have an open embedding
$m\colon D(p_1,\varepsilon_1)^{\beta_1}\times D(p_2,\varepsilon_2)^{\beta_2}
\hookrightarrow\BA^{\beta_1+\beta_2}$ and a canonical isomorphism
$\psi\colon m^*\CI^{\beta_1+\beta_2}\iso\CI^{\beta_1}\boxtimes\CI^{\beta_2}$.
We denote by $r$ the closed embedding $\BA_\BR^\beta\hookrightarrow\BA^\beta$;
we keep the same notation for $D(p,\varepsilon)^\beta_\BR\hookrightarrow
D(p,\varepsilon)^\beta$ in case $p\in\BR$.
We consider the {\em real hyperbolic stalk} $\Phi_\BR(\CI^\beta):=
H^\bullet_c(\BA_\BR^\beta,r^*\CI^\beta)$. According to~\cite[Theorem~I.3.9]{bfs},
$\Phi_\BR(\CI^\beta)$ lives in cohomological degree 0. According
to~\cite[Theorem~I.3.5]{bfs}, we have a canonical isomorphism
$\Phi_\BR(\CI^\beta)^*\simeq\Phi_\BR(\CalD\CI^\beta)$ where $\CalD$ stands for
the Verdier duality.
We have canonical isomorphisms $\Phi_\BR(\CI^\beta)\simeq
H^\bullet_c(D(p,\varepsilon)_\BR^\beta,r^*\CI^\beta)$ for arbitary $p\in\BR,\
\varepsilon\in\BR_{>0}$.

The isomorphism $\psi^{-1}$ above gives rise to the multiplication map
$\Phi_\BR(\CI^{\beta_1})\otimes\Phi_\BR(\CI^{\beta_2})\simeq
H^\bullet_c(D(1,\varepsilon_1)_\BR^{\beta_1},r^*\CI^{\beta_1})\otimes
H^\bullet_c(D(0,\varepsilon_2)_\BR^{\beta_2},r^*\CI^{\beta_2})\to
H^\bullet_c(\BA_\BR^{\beta_1+\beta_2},r^*\CI^{\beta_1+\beta_2})=\Phi_\BR(\CI^{\beta_1+\beta_2})$.
The above selfduality gives rise to the comultiplication map
$\Phi_\BR(\CI^{\beta_1+\beta_2})\to\Phi_\BR(\CI^{\beta_1})\otimes\Phi_\BR(\CI^{\beta_2})$.
According to~\cite[I,II]{bfs}, the twisted graded Hopf algebra
$\Phi_\BR(\CI):=\bigoplus_{\beta\in\BN[I]}\Phi_\BR(\CI^\beta)$ is naturally
isomorphic to $\fu^-$, the negative part of the small quantum group at
$v=\zeta$.

\ssec{rev2}{A review of~\cite{bfs}: factorizable sheaves}
We have an open subset
$\overset{\circ}\BA{}^\beta\supset\overset{\diamondsuit}\BA{}^\beta$ of
configurations of distinct coloured points in $\BA^1\setminus\{0\}$.
It carries a 1-dimensional cohesive local system $\CJ^\beta_\lambda$ with the
monodromies around diagonals same as the ones of $\CJ^\beta$, and also the 
monodromy $\zeta^{2\lambda\cdot\alpha_i}$ when a point of colour $i$ goes around 0
counterclockwise (here $\lambda$ is a weight). We denote by
$\CI^\beta_\lambda$ the Goresky-MacPherson extension of $\CJ^\beta_\lambda$ to
$\BA^\beta$ (a perverse sheaf). Denoting by $A(p,\varepsilon)$ the
complement in $\BA^1$ to the closure of $D(p,\varepsilon)$ (an open annulus),
we have an open embedding $m\colon A(0,\varepsilon)^{\beta_1}\times
D(0,\varepsilon)^{\beta_2}\hookrightarrow\BA^{\beta_1+\beta_2}$ and a canonical
isomorphism $\psi\colon m^*\CI^{\beta_1+\beta_2}_\lambda\iso\CI^{\beta_1}_{\lambda-\beta_2}
\boxtimes\CI^{\beta_2}_\lambda$. A {\em factorizable sheaf} of highest weight
$\lambda$ is a collection of perverse sheaves $\CM^\beta$ on $\BA^\beta$
equipped with factorization isomorphisms
$m^*\CM^{\beta_1+\beta_2}\iso\CI^{\beta_1}_{\lambda-\beta_2}
\boxtimes\CM^{\beta_2}$. In particular, since for $p_1\in\BR$ big enough, and
$\varepsilon_1$ small enough, $D(p_1,\varepsilon_1)\subset A(0,\varepsilon)$,
and the restriction of $\CI^{\beta_1}_{\lambda-\beta_2}$ from
$A(0,\varepsilon)^{\beta_1}$ to
 $D(p_1,\varepsilon_1)^{\beta_1}$ is canonically isomorphic to $\CI^{\beta_1}$, we
obtain isomorphisms
\eq{fact}
\CM^{\beta_1+\beta_2}|_{D(p_1,\varepsilon_1)^{\beta_1}\times
D(0,\varepsilon)^{\beta_2}}\iso\CI^{\beta_1}\boxtimes\CM^{\beta_2}
\end{equation}

Let $a\colon \BA^\beta\to\BA^1$ be the addition, and $\Phi_a(\CM^\beta)$ the
corresponding vanishing cycles. It is a perverse sheaf on the hypersurface
$a=0$, but since $\CM^\beta$ is smooth along coordinate-diagonal stratification,
$\Phi_a(\CM^\beta)$ is supported at the origin $\{\beta\cdot0\}\subset\BA^\beta$,
so we will view $\Phi_a(\CM^\beta)$ just as a vector space. It is canonically
isomorphic to $H^\bullet_c(\BA^\beta_{\BR+},r_+^*\CM^\beta)$ where
$r_+\colon \BA^\beta_{\BR+}\hookrightarrow\BA^\beta$ is the closed embedding of the
``real halfspace'' formed by the real configurations in the preimage
$a^{-1}(\BR_{\geq0})$. Since the vanishing cycles commute with duality, we have a
canonical isomorphism $\Phi_a(\CM^\beta)^*\simeq\Phi_a(\CalD\CM^\beta)$ 
(see e.g.~\cite[Theorem~0.6.3]{bfs}).
The isomorphism~\refe{fact} gives rise to the map
$\Phi_\BR(\CI^{\beta_1})\otimes\Phi_a(\CM^{\beta_2})\simeq
H^\bullet_c(D(p_1,\varepsilon_1)_\BR^{\beta_1},r^*\CI^{\beta_1})\otimes
H^\bullet_c(D(0,\varepsilon)^{\beta_2}_{\BR+},r_+^*\CM^{\beta_2})\to
H^\bullet_c(\BA^{\beta_1+\beta_2}_{\BR+},r_+^*\CM^{\beta_1+\beta_2})=
\Phi_a(\CM^{\beta_1+\beta_2})$, i.e. to the action of $\fu^-\simeq\Phi_\BR(\CI)$ on
$\Phi_a(\CM):=\bigoplus_{\beta\in\BN[I]}\Phi_a(\CM^\beta)$.
The above selfduality gives rise to the coaction
$\Phi_a(\CM^{\beta_1+\beta_2})\to
\Phi_\BR(\CI^{\beta_1})\otimes\Phi_a(\CM^{\beta_2})$; equivalently,
$\Phi_\BR(\CI^{\beta_1})^*\otimes\Phi_a(\CM^{\beta_1+\beta_2})\to\Phi_a(\CM^{\beta_2})$.
Taking into account the isomorphism $\Phi_\BR(\CI^{\beta_1})^*=(\fu^-_{-\beta_1})^*
\simeq\fu^+_{\beta_1}$, we obtain an action of $\fu^+$ on $\Phi_a(\CM)$.
We assign to $\Phi_a(\CM^\beta)$ the weight $\lambda-\beta$. This, together
with the action of $\fu^\pm$, defines the action of $\dot\fu$ on $\Phi_a(\CM)$
(an isomorphism of $\dot\fu$ and the Lusztig small quantum group 
$\dot\bu=\dot\bu_D$ is established in~\cite[Theorem~2.13]{bfs}, 
cf.~\refss{err}).
The resulting functor from the category $\CF\CS$ of factorizable sheaves to
the category $\CalC$ of $\dot\bu$-modules (to be denoted $\Phi$) is an
equivalence of categories.

\ssec{Coxeterfs}{A Coxeter structure on $\CF\CS$}
The diagram $D$ is
the Dynkin diagram of an irreducible Cartan datum $(I,\cdot)$ of finite type.
The category of factorizable sheaves introduced in~\cite[0.4.6]{bfs} will
be denoted by $\CF\CS_D$. For a subdiagram $D'\subset D$ we denote by
$\CF\CS_{D'}$ a similarly defined category with grading by the weight lattice
$X$. That is, compared to the definition of $\CF\CS_D$, the lattice of weights
is always $X$, while the set of colors is $D'\subset D=I$. The braided
balanced tensor structure on $\CF\CS_{D'}$ is introduced
in~\cite[0.5.9,~0.5.10,~IV.6.6]{bfs}.

In order to construct the local systems $F^{\CF\CS}_{D'D''}$ of restriction
functors $\CF\CS_{D'}\to\CF\CS_{D''}$ we
vary the definition~\cite[0.6.7,~0.6.8]{bfs} of the vanishing cycles
functor $\Phi$ in the following way.
Let $\BN[D']\ni\beta=\sum_{j\in D'}b_j\alpha_j$, and let
$\BA^\beta=\prod_{j\in D'}(\BA^1)^{(b_j)}$ be the configuration space of $D'$-colored
effective divisors on the affine line $\BA^1$ with coordinate $t$, of degree
$\nu$. For $D''\subset D'$ we have a pairing
$\langle\cdot,\cdot\rangle\colon \fh_{D'/D''}\times\BA^\beta\to\BA^1$
given in the coordinates $(t_{j,s})_{j\in D',\ 1\leq s\leq b_j}$ on $\BA^\beta$,
and $(z_j)_{j\in D'-D''}$ in the basis of fundamental coweights on $\fh_{D'/D''}$ as follows:
$\langle(z_j),(t_{j,s})\rangle:=\sum_{j\in D'-D''}z_j\sum_{s=1}^{b_j}t_{j,s}$.
The decomposition 
$\beta=\beta''+\ '\beta:=\sum_{j\in D''}b_j\alpha_j+\sum_{k\in D'-D''}b_k\alpha_k$
gives rise to the direct product decomposition
$\BA^\beta=\BA^{\beta''}\times\BA^{'\beta}$. Clearly, given a perverse sheaf 
$\CM^\beta$ on $\BA^\beta$, the vanishing cycles 
$\Phi_{\langle\cdot,\cdot\rangle}\CM^\beta$ is
a perverse sheaf supported on $\fh_{D'/D''}\times\BA^{\beta''}\times0^{'\beta}\simeq
\fh_{D'/D''}\times\BA^{\beta''}$.

\vspace{.3cm}

Let us write $\zeta$ in the form  $\zeta=\exp(\pi i\varkappa)$.

\th{smooth}
If $\CM^\beta$ is a part of data of a
factorizable sheaf $\CM$ and $\varkappa$ is sufficently close
to $0$,  then
$\Phi_{\langle\cdot,\cdot\rangle}\CM^\beta|_{\fh_{D'/D''}^\reg\times\BA^{\beta''}}$ 
is smooth along
$\fh_{D'/D''}^\reg$.
\eth

\prf
It suffices to consider an irreducible $\CM$, and hence an irreducible 
$\CM^\beta$. We may and will assume $D'=D$. Then $\CM^\beta$ is isomorphic to
the Goresky-MacPherson sheaf $\CI^\beta_\lambda$ of~\refss{rev2} for a certain
weight $\lambda$. For $\beta=\sum_{i\in I}b_i\alpha_i$, we consider an 
{\em unfolding} $\pi\colon J\twoheadrightarrow I$ such that for any 
$i\in I,\ \sharp\pi^{-1}(i)=b_i$. Then the product of symmetric groups
$\Sigma_\pi:=\prod_{i\in I}\fS_{b_i}$ acts on the affine space $\BA^J$, and
$\BA^\beta=\BA^J/\Sigma_\pi$. We denote the natural projection 
$\BA^J\to\BA^\beta$ by $\pi$ as well. We denote by 
$\overset{\diamondsuit}\BA{}^J\subset\BA^J$ the complement to the diagonals
in $(\BA^1\setminus\{0\})^J$. We consider the one-dimensional local system 
$\CJ^J_\lambda$ on $\overset{\diamondsuit}\BA{}^J$ with the following monodromies:
$\zeta^{-2\alpha_i\cdot\alpha_j}$ when a point of colour $i$ goes
counterclockwise around a point of colour $j\ne i$, and
$\zeta^{2\lambda\cdot\alpha_i}$ when a point of colour $i$ goes around 0
counterclockwise. We denote by $\CI^J_\lambda$ the Goresky-MacPherson extension
of $\CJ^J_\lambda$ to $\BA^J$. Then $\CI^J_\lambda$ carries an evident
$\Sigma_\pi$-equivariant structure, so $\Sigma_\pi$ acts on the perverse
sheaf $\pi_*\CI^J_\lambda$, and $\CI^\beta_\lambda$ is nothing but the subsheaf
$(\pi_*\CI^J_\lambda)^{\Sigma_\pi,-}$
of $\Sigma_\pi$-antiinvariants in $\pi_*\CI^J_\lambda$ (see~\cite[II.6.13]{bfs}).

For $D''\subset D=I$, let $J'':=\pi^{-1}(D'')\subset J$, and 
$'\!\!J=J\setminus J''$, so that 
$\BA^{J''}\hookrightarrow\BA^J$ (we set the remaining coordinates to be all 
zeros). The construction of~\refss{Coxeterfs} gives rise to a linear map
$m\colon \fh_{D/D''}\times\BA^{J''}\to T^*_{\BA^{J''}}\BA^J=
(\BA^{'\!\!J})^*\times\BA^{J''}$, and we have 
$\Phi_{\langle\cdot,\cdot\rangle}\CM^\beta=
(\pi_*m^\circ\bmu_{\BA^{J''}/\BA^J}\CI^J_\lambda)^{\Sigma_\pi,-}$
where $\bmu$ is the microlocalization functor~\cite{ks}. 

Now by Kashiwara-Schapira theorem identifying the microlocalization and the
Fourier transform 
(see~\cite[Proposition~8.6.3]{ks},~\cite{bryl},~\cite{ks-fourier}), we have
$\bmu_{\BA^{J''}/\BA^J}\CI^J_\lambda\simeq\on{FT}_{\BA^{J''}}\Sp_{\BA^{J''}}\CI^J_\lambda$.
The specialization $\Sp_{\BA^{J''}}\CI^J_\lambda$ is a $\Sigma_\pi$-equivariant
perverse sheaf all of whose simple constituents are of the form 
$\CI^{'\!\!J_1}_{\lambda_1}\boxtimes\CI^{J''_2}_{\lambda_2}$ for certain subsets 
$'\!\!J_1\subset\ '\!\!J,\ J''_2\subset J''$. It suffices to consider
these simple constituents. We have 
$\on{FT}_{\BA^{J''}}(\CI^{'\!\!J_1}_{\lambda_1}\boxtimes\CI^{J''_2}_{\lambda_2})\simeq
(\on{FT}_0\CI^{'\!\!J_1}_{\lambda_1})\boxtimes\CI^{J''_2}_{\lambda_2}$.

The Fourier transform $('\!\pi_*(m^\circ\on{FT}_0
\CI^{'\!\!J_1}_{\lambda_1})|_{\fh_{D/D''}^\reg})^{\Sigma_{'\!\pi},-}$ 
in De Rham setting (i.e.\ the corresponding $D$-module) is calculated 
in~\cite[Theorem~3.2]{fmtv}. It is a version of Casimir connection
\eq{casi}
\nabla=d-\varkappa\sum_{\alpha\in '\!R^+}\frac{\alpha\cdot\alpha}{2}
\frac{d\alpha}{\alpha}f_\alpha e_\alpha
\end{equation}
along $\fh_{D/D''}^\reg$. In particular, the required smoothness assertion 
follows.




In fact, the authors of~\cite{fmtv} work with the De Rham complex of 
``flag logarithmic forms'' but~\cite[Corollary~6.10]{kv} 
implies that if $\varkappa$ is sufficently close to $0$,  this
subcomplex is quasiisomorphic to the De Rham complex $DR(\CI_{\lambda_1}^{'\!\!J})$.
\epr


\conj{smoothconj} The conclusion of \reft{smooth} holds true for
an arbitrary $\zeta\in\BC^*$.
\econj


\ssec{locsys}{}
The collection of perverse sheaves
$\Phi_{\langle\cdot,\cdot\rangle}\CM^\nu|_{\fh_{D'/D''}^\reg\times\BA^{\nu''}},\ \nu\in Y^+$,
enjoys the $D''$-factorization property, i.e. may be viewed as a local
system over $\fh_{D'/D''}^\reg$ of objects of $\CF\CS_{D''}$. This
provides the desired
construction of the local systems $F^{\CF\CS}_{D'D''}$ of restriction functors
$\CF\CS_{D'}\to\CF\CS_{D''}$ required by ~\refd{Coxeter}(Ab),
in the formal neighbourhood of $0$ with respect to $\varkappa$.

If \refco{smoothconj} is true, we get these local systems for every
$\zeta$.

\ssec{iter}{Iterated vanishing cycles.}
The isomorphisms of~\refd{Coxeter}(Ac) are a particular case of the following
construction.
Let $\langle\cdot,\cdot\rangle_W\colon W\times V\to\BA^1$ be a bilinear pairing
between two complex vector spaces. Let $U\subset W$ be a linear subspace.
We denote the restriction of $\langle\cdot,\cdot\rangle$ to $U\times V$
by $\langle\cdot,\cdot\rangle_U$. Let $\CM$ be a perverse sheaf on $V$ smooth
along a
central hyperplane arrangement. We will view $\Phi_{\langle\cdot,\cdot\rangle_W}\CM$
as a perverse sheaf on $W\times W^\perp$ (where $W^\perp\subset V$ is the
annihilator of $W$). Note that the pairing $\langle\cdot,\cdot\rangle_W$
descends to the well defined pairing $\langle\cdot,\cdot\rangle_{W/U}$
between $W/U$ and $U^\perp$.

\th{microsp}
There is a canonical isomorphism
$$
\on{Sp}_{U\times W^\perp}\Phi_{\langle\cdot,\cdot\rangle_W}\CM\iso
\Phi_{\langle\cdot,\cdot\rangle_{W/U}}\Phi_{\langle\cdot,\cdot\rangle_U}\CM
$$
of perverse sheaves on $U\times W/U\times W^\perp$.
\eth

{\bf Proof of~\reft{microsp}.}
Let $V^\perp\subset W$ be the kernel of the bilinear pairing $\langle,\rangle:\
W\times V\to\BA^1$. The smooth base change for the projection $W\to W/V^\perp$
reduces the claim to a construction of a canonical isomorphism
\eq{very old}
\on{Sp}_{Y\times X^\perp/Y\times Y^\perp}\bmu_{Y/V}\CM\iso\bmu_{Y\times X^\perp/X\times X^\perp}
\bmu_{X/V}\CM
\end{equation}
of perverse sheaves on $Y\times(X/Y)^*\times(V/X)^*$,
where $Y:=W^\perp\subset X:=U^\perp\subset V$, and $\bmu$ is the microlocalization functor~\cite{ks}.
This isomorphism will be proved in the next Section, see~\reft{microloc}.

\sec{itera}{Iterated specialization and microlocalization}

\ssec{iterspec}{Iterated specialization}

Fix a complex or real vector space $V$ equipped with a finite  central
hyperplane arrangement $\CH = \{H_i\}$. A linear subspace $V'\subset V$ is called a {\it flat} if it is
an intersection of some hyperplanes from $\CH$. A
filtration
$$
\ldots\subset V_{i+1}\subset V_i\subset\ldots \subset V_0 = V
$$
is called {\it admissible} if all $V_i$ are flats or $0$.

Let $V''\subset V'\subset V$ be an admissible filtration.
$\CH$ induces a central arrangement
on $V'/V''$. Let $Sh_\CH(V'/V'')$ denote the category of constructible sheaves smooth along the corresponding stratification, and let $D^b_\CH(V'/V'')$ denote
the bounded derived category of complexes whose cohomology
belongs to $Sh_\CH(V'/V'')$.

If $V$ is complex and $\CH$ is real, which means that all $H_i$
are given by real equations, then the abelian subcategory
$Perv_\CH(V)\subset D^b_\CH(V)$ admits a description
in terms of linear algebra (quiver)
data, cf.~\cite{ksch}.

Given a flat $W\subset V$, we have the specialization 
functor~\cite{v},~\cite{ks}
$$
\Sp_{W}\colon D^b_\CH(V) \lra D^b_\CH(W\oplus V/W)
$$
whose value on $\CM\in D^b_\CH(V)$ can be described as 
follows~\cite[Section~9]{v}. We fix hermitian metrics on $V$ and on $V/W$.
Let $p$ denote the natural projection $V\to V/W$. 
Let $\xi=(w,u)\in W\oplus V/W$. We fix a sufficiently
small $\varepsilon>0$ and $0<\rho\ll\varepsilon$. We set
$U_{\varepsilon,\rho}=\{v\in V: \rho\|w-v\|+\|\rho u-p(v)\|<\varepsilon\rho\}$:
an open subset of $V$. Then the stalk 
$\Sp_W(\CM)_\xi=R\Gamma(U_{\varepsilon,\rho},\CM)$.


The following lemma is a consequence of this description of 
specialization, cf.~\cite[Theorem~3.17]{fks}.

\lem{third} (a) Let $V_2\subset V_1\subset V$ be an admissible  filtration. 
We have a natural isomorphism of functors $D^b_\CH(V)\lra
D^b_\CH(V_2\oplus V_1/V_2\oplus V/V_1)$,
$$
\phi_{12}\colon \Sp_{V_1/V_2}\Sp_{V_2}\iso
\Sp_{V_2}\Sp_{V_1}.
$$
Let us abbreviate the notation as $\phi_{12}\colon \Sp_{1}\Sp_{2}\iso
\Sp_{2}\Sp_{1}$.

(b) Let $V_3\subset V_2\subset V_1\subset V$ be an admissible  filtration. The various isomorphisms $\phi$ from (a) satisfy the
pentagon relation
$\phi_{23}\circ \phi_{13}\circ \phi_{12} = \phi_{12}\circ \phi_{23}\colon
\Sp_1\Sp_2\Sp_3 \lra \Sp_3\Sp_2\Sp_1.$ \qed
\elem


\sssec{penta}{} Let us explain the relation (b). We have eight categories related to the subquotients of $V$:
$$
D^b_\CH(V)
$$
$$
D^b_\CH(V_3\oplus V/V_3)\ \ \ D^b_\CH(V_2\oplus V/V_2)\ \ \
D^b_\CH(V_1\oplus V/V_1)
$$
$$
D^b_\CH(V_3\oplus V_2/V_3\oplus V/V_2)\ \ \
D^b_\CH(V_3\oplus V_1/V_3\oplus V/V_1)\ \ \
D^b_\CH(V_2\oplus V_1/V_2\oplus V/V_1)
$$
$$
D^b_\CH(V_3\oplus V_2/V_3\oplus V_1/V_2\oplus V/V_1),
$$
they are in bijection with the vertices of a cube. The functors $\Sp$ act 
from a category to all the categories one level below it. There are six 
longest paths from the category $D^b_\CH(V)$ to the category
$D^b_\CH((V_3\oplus V_2/V_3\oplus V_1/V_2\oplus V/V_1)$; 
these paths are in bijection with the symmetric group
$S_3$ (this is an instance of a wellknown geometric fact: the longest paths on 
an $n$-cube are in bijection with the symmetric group $S_n$).

The paths are connected by homotopies arising from
the natural transformations $\phi$; this is a weak Bruhat order on
$S_3$.
Among these six paths there are two which are equal.
The pentagon (b) above  is a hexagon where one of the natural transformations
on the right is the identity.

\ssec{}{Iterated microlocalization} We also have the microlocalization functor 
which may be defined as the composition
$$
\bmu_{W/V}\colon D^b_\CH(V) \overset{\Sp_{W}}\lra D^b_\CH(W\oplus V/W)
\overset{\on{FT}_W}\lra D^b_\CH(W\oplus (V/W)^*)
$$
where $\on{FT}_W$ is the Fourier-Sato transformation~\cite{ks}.
In order to keep track of the ambient space, from now on we will use another
notation for the specialization: $\Sp_W=\Sp_{W/V}$.

\th{microloc} Let $V$ be a complex vector space and let
$Y\subset X\subset V$ be an admissible filtration.
For $\CM\in Perv_\CH(V)$ there exists a canonical isomorphism
\eq{old}
\psi_{XY}\colon \on{Sp}_{Y\times X^\perp/Y\times Y^\perp}\bmu_{Y/V}\CM\iso\bmu_{Y\times X^\perp/X\times X^\perp}
\bmu_{X/V}\CM
\end{equation}
of perverse sheaves on $Y\times(X/Y)^*\times(V/X)^*$.

These isomorphisms satisfy the pentagon relation connected with
an admissible filtration $Y\subset X\subset Z\subset V$.
\eth

\ssec{proo}{Proof of~\reft{microloc}}
The  rest of this Section is devoted to the proof of this theorem.

The following properties are consequences of the definition.

\sssec{commute}{} For a linear subspace $Y\subset V$ and a perverse sheaf
$\CM$ on $V$ we have a canonical isomorphism
$\bmu_{Y/V}\CM\iso\bmu_{Y/Y\times(V/Y)}\on{Sp}_{Y/V}\CM$.

\sssec{fourier}{} For a perverse sheaf $\CM$ on the product of two vector spaces $Y\times Z$, monodromic along the projection $Y\times Z\to Y$,
we have a canonical isomorphism $\bmu_{Y/Y\times Z}\CM\iso FS_{Y/Y\times Z}\CM$
between the microlocalization and the Fourier-Sato transform on the vector bundle $Y\times Z\to Y$.

\sssec{dmod}{}
We will have to work on the $D$-module side of the Riemann-Hilbert
correspondence, so we recall the definition of the specialization functor
in this context. Given a linear subspace $Y\subset V$ we choose a complementary
subspace $Z\subset V,\ V\simeq Y\oplus Z$, with linear coordinates
$z_1,\ldots,z_d$. Then the ring of differential operators $D_V$ has a grading
such that $\deg z_i=1,\ \deg\partial_{z_i}=-1,\ \deg y=\deg\partial_y=0$ for any
$i=1,\ldots,d$, and any linear coordinates $y$ on $Y$. Let $F^\bullet D_V$ be
the corresponding descending filtration. Note that $\on{gr}^FD_V$ is canonically
isomorphic to $D_{Y\times(V/Y)}$.
Let $M$ be a regular holonomic $D$-module on $V$. It possesses a unique
(descending) Malgrange-Kashiwara filtration
$\ldots F^{-1}M\supset F^0M\supset F^1M\supset\ldots$ compatible with the
filtration on $D_V$ such that (a) for $j>0$ and $k\gg0$, we have
$F^{\pm k\pm j}M=F^{\pm j}D_VF^{\pm k}M$; (b) the generalized eigenvalues of the Euler
vector field $\sum_{1\leq i\leq d}z_i\partial_{z_i}$ on $F^kM/F^{k+1}M$ have real
parts in $[k,k+1)$ for any $k\in\BZ$.

The specialization of $M$ is defined as a $D_{Y\times(V/Y)}$-module
$\on{Sp}_{Y/V}M:=\on{gr}^FM$. We say that $M$ is {\em potentially monodromic}
along $Y$ if the Malgrange-Kashiwara filtration $F^\bullet M$ is compatible
with some grading $G^\bullet M$ compatible with the grading on $D_V$.
Equivalently, we can choose a complementary subspace $Z\subset V$ such that
$M$ is monodromic along the corresponding projection $V\to Y$. Any such choice
defines the isomorphisms $Y\times Z\iso V$ and $\on{Sp}_{Y/V}M\iso M$.

\sssec{foumic}{} For a regular holonomic $D$-module on the product of two
vector spaces $Y\times Z$, monodromic along the projection $Y\times Z\to Y$,
its microlocalization $\bmu_{Y/Y\times Z}M\iso FT_{Y/Y\times Z}$ is the following
regular holonomic $D_{Y\times Z^*}$-module. We choose some linear coordinates
$y_1,\ldots,y_e$ on $Y$, and $z_1,\ldots,z_d$ on $Z$. Let $\xi_1,\ldots,\xi_d$
be the dual coordinates on $Z^*$. Then we have an isomorphism
$D_{Y\times Z^*}\iso D_{Y\times Z}\colon y_i\mapsto y_i,\
\partial_{y_i}\mapsto\partial_{y_i},\ \xi_j\mapsto\partial_{z_j},\
\partial_{\xi_j}\mapsto-z_j$ (independent of the choice of coordinates),
and $\bmu_{Y/Y\times Z}M$ is nothing but $M$ viewed
as a $D_{Y\times Z^*}$-module via this isomorphism.

Now we proceed to the construction of isomorphism~\refe{old} in the
$D$-module setting. All the $D_V$-modules below are assumed to be regular
holonomic and smooth along some central hyperplane arrangement in $V$.

\lem{first}
If a regular holonomic $D_V$-module $M$ is potentially monodromic along
both $Y$ and $X$, then there is a canonical isomorphism~\refe{old}.
\elem

\prf The microlocalization $\bmu_{Y/V}M$ is monodromic along the projection
$Y\times Y^\perp\to Y$ and potentially monodromic along $Y\times X^\perp$.
Any choice of the direct complement $Z\subset X\colon X=Y\oplus Z$ as in the
definition of potential monodromicity, gives rise to the isomorphisms
$Y\times Z^*\times X^\perp\iso Y\times Y^\perp$ and
$\on{Sp}_{Y\times X^\perp/Y\times Y^\perp}\bmu_{Y/V}M\iso\bmu_{Y/V}M$.
So it remains to construct an isomorphism
$\bmu_{Y/V}M\iso\bmu_{Y\times X^\perp/X\times X^\perp}\bmu_{X/V}M$. We choose a
direct complement $S\subset V\colon V=X\oplus S$ as in the definition of
potential monodromicity. Then the required isomorphism follows from the
explicit formulas in~\refsss{foumic}. One can check that it does not depend
on the choices of the complements $Z$ and $S$.
\epr

\lem{second}
If a regular holonomic $D_V$-module $M$ is potentially monodromic along
$Y$, then there is a canonical isomorphism~\refe{old}.
\elem

\prf By~\refsss{commute} we can replace the RHS of~\refe{old} by
$\bmu_{Y\times X^\perp/X\times X^\perp}\bmu_{X/V}\on{Sp}_{X/V}M$. But $\on{Sp}_{X/V}M$
is potentially monodromic along both $Y$ and $X$. So it remains to
use~\refl{first} and construct an isomorphism
\eq{93}
\on{Sp}_{Y\times X^\perp/Y\times Y^\perp}\bmu_{Y/V}M\iso\bmu_{Y/X\times(V/X)}\on{Sp}_{X/V}M
\end{equation}
We choose a direct complement $Z\subset X\colon X=Y\oplus Z$, and a direct
complement $S\subset V\colon V=X\oplus S$ such that $Z\oplus S$ is as in the
definition of potential monodromicity. We choose the linear coordinates
$z_1,\ldots,z_d$ in $Z$, and $s_1,\ldots,s_e$ in $S$, and the dual
coordinates $\xi_1,\ldots,\xi_d$ in $Z^*$, and $\eta_1,\ldots,\eta_e$ in $S^*$.
The $D_V$-module $M$ has a grading $G^\bullet M$ compatible with the grading
on $D_V$ such that $\deg z_i=\deg s_j=1,\
\deg\partial_{z_i}=\deg\partial_{s_j}=-1,\ \deg y_k=\deg\partial_{y_k}=0$.
According to~\refsss{foumic}, the microlocalization $\bmu_{Y/V}M$ amounts to the
substitution $\xi_i\mapsto\partial_{z_i},\ \eta_j\mapsto\partial_{s_j},\
\partial_{\xi_i}\mapsto-z_i,\ \partial_{\eta_j}\mapsto-s_j,\ y_k\mapsto y_k,\
\partial_{y_k}\mapsto\partial_{y_k}$.
To compute the specialization
$\on{Sp}_{X/V}M$ we use the Malgrange-Kashiwara filtration $F^\bullet M$ as
in~\refsss{dmod}. To compute the LHS of~\refe{93} we use the unique filtration
$'F^\bullet M$ compatible as in~\refsss{dmod} with the grading on $D_V$ such
that $\deg z_i=-1,\ \deg\partial_{z_i}=1,\ \deg s_j=\deg\partial_{s_j}=\deg y_k=
\deg\partial_{y_k}=0$ (note that it is {\em not} the Malgrange-Kashiwara
grading of $D_V$, but rather the Fourier image of one).

The construction of the isomorphism~\refe{93} amounts to the construction
of the isomorphism
\eq{94}
\on{gr}^FM\iso\on{gr}^{'F}M
\end{equation}
Note that both $F^\bullet M$ and $'F^\bullet M$ are compatible with the grading
$G^\bullet M$, that is $F^iM=\oplus_j(F^iM\cap G^jM)$ and
$'F^iM=\oplus_j(\ 'F^iM\cap G^jM)$ for any $i\in\BZ$. From the uniqueness of
$'F^\bullet M$ and $F^\bullet M$ it follows that $'F^kM\cap G^jM=F^{k+j}M\cap G^jM$.
Thus the desired isomorphism~\refe{94} is the direct sum of natural
isomorphisms $(\ 'F^kM\cap G^jM)/(\ 'F^{k+1}M\cap G^jM)=
(F^{k+j}M\cap G^jM)/(F^{k+j+1}M\cap G^jM)$.
\epr

\ssec{finconstr}{The end of the proof}
Now we can finish the construction of isomorphism~\refe{old} for an arbitrary
$D_V$-module $M$. It suffices to construct the isomorphism~\refe{93} for
arbitrary $M$ (not necessarily potentially monodromic along $Y$).
By~\refsss{commute} we can replace the RHS of~\refe{93} by
$\bmu_{Y/X\times(V/X)}\on{Sp}_{Y/X\times(V/X)}\on{Sp}_{X/V}M$. By~\refl{third} we
can replace this by
$\bmu_{Y/X\times(V/X)}\on{Sp}_{Y\times(X/Y)/Y\times(V/Y)}\on{Sp}_{Y/V}M$.
However, $\on{Sp}_{Y/V}M$ is already potentially monodromic along $Y$, so
by~\refl{second} we have
$\on{Sp}_{Y\times X^\perp/Y\times Y^\perp}\bmu_{Y/V}\on{Sp}_{Y/V}M\iso
\bmu_{Y/X\times(V/X)}\on{Sp}_{Y\times(X/Y)/Y\times(V/Y)}\on{Sp}_{Y/V}M$. One last
application of~\refsss{commute} allows to replace the LHS by
$\on{Sp}_{Y\times X^\perp/Y\times Y^\perp}\bmu_{Y/V}M$.

This completes  the proof of~\reft{microloc}
and hence that of~\reft{microsp}.

\rem{discrepancy} The statements of~\refl{third} and~\reft{microloc} look
similar, but the proofs are very different: one is topological, another is
De Rham (via $D$-modules). Let us comment on this discrepancy.
On the one hand,~\refl{third} has an easy proof in De Rham setting as well.
On the other hand, 
let $Perv_\CH(V)\subset D_\CH^b(V)$ denote the subcategory of perverse sheaves.
If all $H_i\in \CH$ are given by real equations then $Perv_\CH(V)$
admits an explicit description in terms of linear algebra (quiver)
data, cf.~\cite{ksch}. The specialization and microlocalization functors can be
described in this language, and a topological proof 
of~\refl{third} (resp.~\reft{microloc}) is given in~\cite[Theorem~3.17]{fks} 
(resp.~\cite[Theorem~5.6]{fks}).
\erem

\sec{discu}{Discussion}

\ssec{goals}{Desiderata}
We have gone to all this trouble just to conjecture that the functor $\Phi$
of~\cite{bfs} takes the Coxeter structure of~\refss{Coxeterfs} to the one of~\refss{Coxeterc}.
By~\cite[Theorem~3.2]{fmtv} and Fourier=microlocalization this would imply
that the monodromy of the Casimir connection is given by the Lusztig
symmetries for any $\zeta$. Note that the functor $\Phi$ of~\cite{bfs} is
nothing but the stalk at the fundamental chamber $C^{D'}_0$
of the local system of the restriction functors
$F^{\CF\CS}_{D'\emptyset}\colon \CF\CS_{D'}\to\on{Vect}_X$ of~\refss{Coxeterfs}.

\ssec{ref}{Tilted functors $\Phi$}

Recall the setup of~\refss{rev2}. Let us choose a point $z^{(w)}$ in a chamber
$C_w\subset\fh^\reg_\BR,\ w\in W$. Instead of $\Phi_a(\CM^\beta)$
let us consider the spaces of ``tilted"
vanishing cycles $\Phi_w(\CM^\beta):=\Phi_{\langle z^{(w)},?\rangle}(\CM^\beta)\simeq
H^\bullet_c(\BA^\beta_{\BR w},r_w^*\CM^\beta)$ where
$r_w\colon \BA^\beta_{\BR w}\hookrightarrow\BA^\beta$ is the closed embedding of the
``tilted halfspace'' formed by the real configurations in the preimage
$\langle z^{(w)},?\rangle^{-1}(\BR_{\geq0})$. In particular, $\Phi_a\simeq\Phi_e$.

For a real $p_1$ such that
$\langle z^{(w)},\beta\cdot p_1\rangle$ is positive and big enough, similarly
to~\refss{rev2} we obtain the map
$$
\Phi_\BR(\CI^{\beta_1})\otimes\Phi_w(\CM^{\beta_2})\simeq
H^\bullet_c(D(p_1,\varepsilon_1)_\BR^{\beta_1},r^*\CI^{\beta_1})\otimes
H^\bullet_c(D(0,\varepsilon)^{\beta_2}_{\BR w},r_w^*\CM^{\beta_2})
\to
$$
$$
H^\bullet_c(\BA^{\beta_1+\beta_2}_{\BR w},r_w^*\CM^{\beta_1+\beta_2})=
\Phi_w(\CM^{\beta_1+\beta_2}),
$$
i.e. the action of $\fu^-\simeq\Phi_\BR(\CI)$ on
$\Phi_w(\CM):=\bigoplus_{\beta\in\BN[I]}\Phi_w(\CM^\beta)$.

The selfduality of $\Phi_w$ gives rise to the action of $\fu^+$ on
$\Phi_w(\CM)$ similarly to~\refss{rev2}.
We assign to $\Phi_w(\CM^\beta)$ the weight $w(\lambda-\beta)$, and using
the isomorphisms $T'_{w\pm}\colon \fu^+\iso T'_{w\pm}(\fu^+)\subset\fu,\
\fu^-\iso T'_{w\pm}(\fu^-)\subset\fu$ we obtain the action of $T'_{w\pm}(\fu^+),
T'_{w\pm}(\fu^-)$ on $\Phi_w(\CM)$. This, together with the above grading,
defines an action of $\dot\fu$ on $\Phi_w(\CM)$, i.e. gives rise to two functors $\Phi_{w\pm}\colon \CF\CS\to\CalC$.

\bigskip

Given a straight line interval $\gamma_w$ from
$z^{(e)}$ to $z^{(w)}$ we obtain the corresponding ``halfmonodromy''
transformations $\gamma_{w,\beta}^\pm(\CM)\colon \Phi_e(\CM^\beta)\iso\Phi_{w\pm}(\CM^\beta)$
(independent of the choice of $\gamma$) for any
factorizable sheaf $\CM$ and $\beta\in\BN[I]$.

The following conjecture is a reformulation of~\refss{goals}.

\conj{refor} The maps $\{\gamma_{w,\beta}^\pm(\CM),\ \CM\in \CF\CS\}$ define two  natural transformations of functors
$\gamma_w^+\colon \Phi\iso\Phi_{w +}$ and $\gamma_w^-\colon \Phi\iso\Phi_{w -}$.
\econj

As we have already mentioned, the main theorem of~\cite{tr} implies this 
conjecture for $\zeta=\exp(\hbar),\ \hbar$ being a formal parameter.

\bigskip

\footnotesize{ {\bf M.F.}: National Research University Higher
School of Economics, Russian Federation,\\ 
Math. Dept., 6 Usacheva st, Moscow 119048;\\
Skolkovo Institute of Science and Technology;\\
Institute for Information Transmission Problems of RAS\\ 
{\tt fnklberg@gmail.com}}

\footnotesize{ {\bf V.S.}: Institut de Math\'ematiques de Toulouse,
Universit\'e Paul Sabatier, 31062 Toulouse, France;\\ 
{\tt schechtman@math.ups-tlse.fr}}

\end{document}